\documentclass[aop,MSNbibl,seceqn,dvips]{arximspdf}

\doi{10.1214/12-AOP830} %kopijuoti is PTS
\volume{42}
\issue{2}
\pubyear{2014}
\firstpage{576}
\lastpage{622}

\makeatletter

\def\cal{\mathcal}

\newtheorem{theorem}{Theorem}[section]
\newtheorem{corollary}[theorem]{Corollary}
\newtheorem{lemma}[theorem]{Lemma}
\newtheorem{lemmas}{Lemma}[section]
\newtheorem{proposition}[theorem]{Proposition}
\newproclaim{definition}[theorem]{Definition}
\newproclaim{assumption}[theorem]{Assumption}
\newproclaim{remark}[theorem]{Remark}

\def\R{\mathbb{R}}
\def\Z{\mathbb{Z}}
\def\E{\mathbb{E}}
\def\P{\mathbb{P}}

\def\ol{\overline}
\def\Var{\operatorname{Var}}
\def\Cov{\operatorname{Cov}}

\makeatother

\begin{document}
\begin{frontmatter}

\title{Quenched asymptotics for Brownian motion in
generalized Gaussian potential}
\runtitle{Brownian motion in generalized Gaussian potential}

\begin{aug}
\author[A]{\fnms{Xia} \snm{Chen}\corref{}\thanksref{t1}\ead[label=e1]{xchen@math.utk.edu}}
\thankstext{t1}{Supported in part by the Simons Foundation \#244767.}
\runauthor{X. Chen}
\affiliation{University of Tennessee}
\address[A]{Department of Mathematics\\
University of Tennessee\\
Knoxville, Tennessee 37996\\
USA\\
\printead{e1}} %adresu isvedimo komanda gale!
\end{aug}

% HISTORY:
\received{\smonth{1} \syear{2012}}
\revised{\smonth{8} \syear{2012}}

% ABSTRACT
%
\begin{abstract}
In this paper,
we study the long-term asymptotics for the quenched moment
\[
\E_x\exp \biggl\{\int_0^tV(B_s)\,ds
\biggr\}
\]
consisting of a $d$-dimensional Brownian motion $\{B_s; s\ge 0\}$
and a generalized Gaussian field $V$. The major progress made
in this paper includes: Solution to
an open problem posted by Carmona and Molchanov
[\textit{Probab. Theory Related Fields} \textbf{102} (1995)
433--453], the quenched
laws for Brownian motions in Newtonian-type potentials and in the potentials
driven by white noise or by fractional white noise.
\end{abstract}

% KEYWORDS

\end{frontmatter}

%s1 #&#
\section{Introduction} \label{intro}

The classic Anderson model can be formulated
as the following heat equation:
%
%
%e1.1 #&#
%
\begin{equation}
\label{intro-1} \cases{ %
\partial_tu(t,x)
=\frac{1}{2}\Delta u(t,x)+V(x)u(t,x),
\vspace*{2pt}\cr
u(0, x)=1,}
\end{equation}
where $\{V(x);  x\in\R^d\}$ is often made
as a stationary
random field called potential.

Under some regularity assumption such as H\"older continuity on $V(x)$,
the system has a unique solution with Feynman--Kac representation
%
%
%e1.2 #&#
%
\begin{equation}
\label{intro-2} u(t,x)=\E_x\exp \biggl\{\int_0^{t}V(B_{s})\,ds
\biggr\},
\end{equation}
where $\{B_t; t\ge0\}$ is
a $d$-dimensional Brownian motion independent of
$V(x)$, and~$\E_x$ is the expectation with respect to
$B_t$ given $B_0=x$.

An important aspect in studying parabolic Anderson models is its
long-term asymptotics. There are two types of asymptotics: one is
labeled as quenched law concerning the limit behavior of
the random field $u(t,x)$ conditioning on the random potential
$V(x)$; another is known as annealed law
with interest in the limit behavior of
$\E u(t,x)$ and other deterministic moments
of $u(t,x)$.
In the case when $\{V(x);  x\in\R^d\}$
is a mean zero
stationary Gaussian field with the covariance function
%
%
%e1.3 #&#
%
\begin{equation}
\label{intro-5} \gamma(x)=\Cov \bigl(V(0), V(x) \bigr),\qquad x\in\R^d.
\end{equation}
Carmona and Molchanov
(Theorem~5.1, \cite{CM-1}) establish the quenched law
%
%
%e1.4 #&#
%
\begin{equation}
\label{intro-6} \lim_{t\to\infty}\frac{1}{ t\sqrt{\log t}} \log
\E_x\exp \biggl\{\int_0^tV(B_s)\,ds
\biggr\} =\sqrt{2d\gamma(0)},\qquad \mbox{a.s.}
\end{equation}
under the condition
$ \lim_{\vert x\vert\to\infty}\gamma(x)=0$.
See \cite{G-K-M} for the asymptotics of the second order,
and
\cite{BK1,CV,G-M-1,Sznitman}
and \cite{VZ}
for a variety of versions in literature.

This paper is concerned with the setting of the generalized Gaussian
fields, in which the potential $V$ is not defined pointwise.
A typical example is when $V$ is a white or fractional
white noise.
Recall that a generalized function
is defined as a linear functional
$\{\langle\xi, \varphi\rangle;  \varphi\in
{\cal S}(\R^d)\}$ on a suitable
space ${\cal S}(\R^d)$
of the functions known as the test functions.
The classic notion of function is generalized in the sense that
%
%
%e1.5 #&#
%
\begin{equation}
\label{intro-12} \langle\xi,\varphi\rangle=\int_{\R^d}\xi(x)
\varphi(x)\,dx ,\qquad\varphi\in{\cal S}\bigl(\R^d\bigr)
\end{equation}
whenever $\xi(x)$ is a ``good'' function defined pointwise on $\R^d$.
We refer the book~\cite{GV} by Gel'fand and Vilenkin
for details.

A generalized random field $V$ is a generalized random function.
In this paper, we consider the case when
${\cal S}(\R^d)$ is the
Schwartz
space of rapidly
decreasing and infinitely smooth functions, and
$\{\langle V,\varphi\rangle; \varphi\in{\cal S}(\R^d)\}$
is a mean-zero Gaussian field satisfying the homogeneity
%
%
%e1.6 #&#
%
\begin{equation}
\label{intro-9} \bigl\{\bigl\langle V, \varphi(\cdot-x)\bigr\rangle; \varphi\in{
\cal S}\bigl(\R^d\bigr)\bigr\}\stackrel {d} { =}\bigl\{\langle V,
\varphi\rangle; \varphi\in{\cal S}\bigl(\R^d\bigr)\bigr\},\qquad x\in
\R^d.
\end{equation}
The covariance functionals
$\Cov(\langle V,\varphi\rangle, \langle V, \psi\rangle)$
of the generalized Gaussian fields considered in this work are
continuous on ${\cal S}(\R^d)\times
{\cal S}(\R^d)$. Consequently, $\{\langle V,\varphi\rangle;
\varphi\in{\cal S}(\R^d)\}$ is
continuous in probability and therefore yields a measurable version.

The classic Bochner representation can be generalized
((1), page 290, \cite{GV}) in the following way: There
is a positive measure $\mu(d\lambda)$ on $\R^d$,
known as spectral measure, such that
%
%
%e1.7 #&#
%
\begin{equation}
\label{intro-25} \Cov\bigl(\langle V,\varphi\rangle, \langle V, \psi\rangle\bigr)
=\frac{1}{(2\pi)^d}\int_{\R^d}{\cal F}(\varphi) (\lambda)\ol{{
\cal F}(\psi) (\lambda)}\mu(d\lambda),
\end{equation}
where ${\cal F}(\varphi)(\lambda)$ denotes the Fourier transform
of the function $\varphi\in{\cal S}(\R^d)$.
Further, $\mu(d\lambda)$ is tempered in the sense that
$(1+ \vert\cdot\vert^2)^{-p}\in{\cal L}(\R^d, \mu)$
for some $p>0$.

In the settings considered in this paper, the notion of covariance function
$\gamma(\cdot)$ defined by (\ref{intro-5}) can also be extended to
the form
%
%
%e1.8 #&#
%
\begin{equation}
\label{intro-10}\qquad \Cov\bigl(\langle V,\varphi\rangle, \langle V, \psi\rangle\bigr)
=\int_{\R^d}\gamma(x-y)\varphi(x)\psi(y) \,dx\,dy,\qquad \varphi, \psi
\in{\cal S}\bigl(\R^d\bigr)
\end{equation}
with $\gamma(x)=\delta_0(x)$ (Dirac function)
or with $\gamma(x)$ being defined pointwise
on $\R^d\setminus\{0\}$ and satisfying $\gamma(0)\equiv
\lim_{x\to0}\gamma(x)=\infty$---in both cases
$\mu(d\lambda)$ is an infinite measure.
As a consequence, it is impossible to make $V$ pointwise defined through
relation (\ref{intro-12}),
for otherwise we would have to face the ``Gaussian variable'' $V(x)$
with $\Var (V(x) )=\gamma(0)=\infty$ for every $x\in\R^d$.

Nevertheless, representation (\ref{intro-2}) can be extended to
the generalized setting under some suitable condition. The generalized
Gaussian potentials appearing in our main theorems satisfy
(Lemma~\ref{intro-25'})
%
%
%e1.9 #&#
%
\begin{equation}
\label{intro-50} \int_{\R^d}\frac{1}{(1+ \vert\lambda\vert^2)^{1-\delta}}\mu (d\lambda )<
\infty
\end{equation}
for some $\delta>0$. As a consequence
(Lemma~\ref{def}), the $L^2$-limit
\[
\int_0^t V (B_s)\,ds\stackrel{
\mathrm{def}} { =}\lim_{\varepsilon\to0^+}\int_0^t
V_\varepsilon(B_s)\,ds
\]
exists and, the time integral defined in this way yields a continuous
version as a stochastic process, where the pointwise
defined Gaussian field $V_\varepsilon(x)$
appears as a smoothed version of
$V$; see Lemma~\ref{def} for details. In
addition, the time integral defined in this way
is exponentially integrable with respect to
$\E_x$, as pointed out in Section~\ref{u}.
Consequently, representation
(\ref{intro-2}) makes sense in our settings.
According to a treatment proposed on page 448 of \cite{CM-1},
it solves the Anderson model (\ref{intro-1}) in some proper sense.
The major goal
of this work is to study the large-$t$ behavior of the quenched
exponential moment in (\ref{intro-2}).

In \cite{CM-1},
Carmona and Molchanov ask what happens when
the covariance function $\gamma(x)$ is defined pointwise, continuous
in $\R^d\setminus\{0\}$ but $\gamma(0)=\infty$ with the degree of
singularity
measured by
%
%
%e1.10 #&#
%
\begin{equation}
\label{intro-15} \gamma(x)\sim c(\gamma)\vert x\vert^{-\alpha}\qquad (x\to0)
\end{equation}
for some $0<\alpha<2$ and $c(\gamma)>0$.
Here we point out that the restriction ``$\alpha<d$''
has to be added for the covariance functional
$\Cov (\langle V, \varphi\rangle, \langle V, \psi\rangle )$
to be well-defined. Indeed,
for a nonnegative $\varphi\in{\cal S}(\R^d)$
strictly positive in a neighborhood of $0$, there are $C>0$ and
$\varepsilon>0$
such that
\[
\Var \bigl(\langle V, \varphi\rangle \bigr)\ge C^{-1} \int
_{\{\vert x\vert\le\varepsilon\}\times\{\vert y\vert\le\varepsilon
\}} \frac{dx\,dy}{\vert x-y\vert^\alpha}.
\]
The right-hand side diverges if $\alpha\ge d$.\vadjust{\goodbreak}

In their paper, Carmona and Molchanov \cite{CM-1}
conjecture that under (\ref{intro-15}),
\[
\log\E_x\exp \biggl\{\int_0^tV(B_s)\,ds
\biggr\} \approx t(\log t)^{{(4-\alpha)}/ {(2-\alpha)}},\qquad \mbox{a.s. } (t\to
\infty).
\]

The following theorem tells a slightly different story.

%
%th1.1 #&#
\begin{theorem}\label{intro-33} Let the covariance function
$\gamma(x)$ be continuous on $\R^d\setminus\{0\}$
and be bounded outside every neighborhood of 0. Assume
(\ref{intro-15}) with $0<\alpha<2\wedge d$. Then for any $x\in\R^d$,
%
%
%e1.11 #&#
%
\begin{eqnarray}
\label{intro-35} &&\lim_{t\to\infty}t^{-1}(\log
t)^{-{2}/{(4-\alpha)}}\log\E_x\exp \biggl\{ \int_0^tV(B_s)\,ds
\biggr\}
\nonumber
\\[-8pt]
\\[-8pt]
\nonumber
&&\qquad=\frac{4-\alpha}{4} \biggl(\frac{\alpha}{2} \biggr)^{{\alpha}/ {(4-
\alpha)}} \bigl(2dc(\gamma)\kappa(d,\alpha) \bigr)^{{2}/ {(4-\alpha)}},\qquad
\mbox{a.s.},
\end{eqnarray}
where the constant $c(\gamma)>0$ is given in (\ref{intro-15}), and
$\kappa(d,\alpha)>0$ is the best constant of the inequality
[see (\ref{a-2}) in the \hyperref[app]{Appendix}]
\[
\int\!\! \int_{\R^d\times\R^d}\frac{f^2(x)f^2(y)}{\vert x-y\vert
^\alpha} \le C\|f
\|_2^{4-\alpha}\|\nabla f\|_2^\alpha,\qquad f\in
W^{1,2}\bigl(\R^d\bigr)
\]
with $W^{1,2}(\R^d)$ being defined as the Sobolev space
%
%
%e1.12 #&#
%
\begin{equation}
\label{intro-31'} W^{1,2}\bigl(\R^d\bigr)=
\bigl\{f\in{\cal L}^2\bigl(\R^d\bigr); \nabla f\in{\cal
L}^2\bigl(\R^d\bigr)\bigr\}.
\end{equation}
\end{theorem}

We now consider a special case. In light
of some classical laws of physics,
such as Newton's gravity law and Coulomb's electrostatics law, it makes
sense to consider the potential formally given as
\[
V(x)=\int_{\R^d}\frac{1}{\vert x-y\vert^p}W(dy), \qquad x\in\R^d
\]
in the parabolic Anderson model (\ref{intro-1}). Here
$\{W(x);  x\in\R^d\}$ is a standard Brownian sheet.
The relevant Gaussian field
%
%
%e1.13 #&#
%
\begin{equation}
\label{intro-17} \langle V, \varphi\rangle =\int_{\R^d}
\biggl[\int_{\R^d}\frac{\varphi(y)}{\vert y-x\vert
^p}\,dy \biggr] W(dx),\qquad \varphi
\in{\cal S}\bigl(\R^d\bigr)
\end{equation}
is well defined with the covariance function
%
%
%e1.14 #&#
%
\begin{equation}
\label{intro-18} \gamma(x)=C(d,p)\vert x\vert^{-(2p-d)},
\end{equation}
provided $d/2<p<\frac{d+2}{2}\wedge d$, where
%
%
%e1.15 #&#
%
\begin{equation}
\label{intro-19} C(d,p)=\pi^{d/2}\frac{ \Gamma^{2} ({(d-p)}/{2} )
\Gamma ({(2p-d)}/{2} )}{
\Gamma^{2} ({p}/{2} )
\Gamma({d-p})}.
\end{equation}
Indeed,
\begin{eqnarray*} \Cov \bigl(\langle V, \varphi\rangle, \langle V,
\psi\rangle \bigr)&=&\int_{\R^d} \biggl[\int_{\R^d}
\frac{\varphi(y)\,dy}{\vert y-x\vert^p} \biggr] \biggl[\int_{\R^d}\frac{\psi(z)\,dy}{\vert z-x\vert^p}
\biggr]\,dx
\\
&=&\int_{\R^d\times\R^d}\varphi(y)\psi(z) \biggl[\int
_{\R^d} \frac{dx}{\vert y-x\vert^p\vert z-x\vert^p} \biggr]\,dy\,dz\\
& =&C(d,p)\int
_{\R^d\times\R^d}\frac{\varphi(y)\psi(z)}{\vert
y-z\vert^{2p-d}}\,dy\,dz,
\end{eqnarray*}
where the last step follows from
(1.31) in \cite{C-R} (with $\sigma$ being replaced by $2p-d$).

Thus (\ref{intro-15}) holds with $\alpha=2p-d<2\wedge d$.

%
%co1.2 #&#
\begin{corollary}\label{intro-21}
In the special case given in (\ref{intro-17}) with $d/2<p<d\wedge
\frac{d+2}{2}$,
%
%
%e1.16 #&#
%
\begin{eqnarray}
\label{intro-22} &&\lim_{t\to\infty}t^{-1}(\log
t)^{-{2}/{(4+d-2p)}}\log\E_x\exp \biggl\{ \theta\int
_0^tV(B_s)\,ds \biggr\}
\nonumber\\
&&\qquad=\frac{4+d-2p}{4} \biggl(\frac{2p-d}{2} \biggr)^{{(2p-d)}/
{(4+d-2p)} }\\
&&\qquad\quad{}\times\bigl(2dC(d,p)\theta^2\kappa(d,2p-d) \bigr)^{{2}/
{(4+d-2p)}}, \qquad\mbox{a.s.}
\nonumber
\end{eqnarray}
for any $\theta>0$, where $C(d,p)>0$ is given in (\ref{intro-19}).
\end{corollary}

In the next theorem, the potential is a fractional
white noise formally written
as
\[
V(x)=\frac{\partial^d W^H}{\partial x_1\cdots\partial x_d}(x),\qquad x=(x_1,\ldots, x_d)\in
\R^d,
\]
where
$W^H(x)$ ($x=(x_1,\ldots, x_d)\in\R^d$) is a fractional Brownian sheet
with Hurst index $H=(H_1,\ldots, H_d)$. We assume that
%
%
%e1.17 #&#
%
\begin{equation}
\label{intro-37} \frac{1}{2}<H_j<1\qquad (j=1,\ldots, d)\quad \mbox
{and}\quad \sum_{j=1}^dH_j>d-1.
\end{equation}
The generalized Gaussian field relevant to the
problem is defined by the stochastic integral
%
%
%e1.18 #&#
%
\begin{equation}
\label{intro-38} \langle V,\varphi\rangle =\int_{\R^d}
\varphi(x)W^H(dx),\qquad \varphi\in{\cal S}\bigl(\R^d\bigr).
\end{equation}
In this setting,
%
%
%e1.19 #&#
%
\begin{eqnarray}
\label{intro-39} \gamma(x)&=&C_H \Biggl(\prod
_{j=1}^d\vert x_j\vert^{2-2H_j}
\Biggr)^{-1}\quad \mbox{and}
\nonumber
\\[-8pt]
\\[-8pt]
\nonumber
 \mu(d\lambda)&=&\widehat{C}_H \Biggl(
\prod_{j=1}^d\vert\lambda_j
\vert^{2H_j-1} \Biggr)^{-1} \,d\lambda,
\end{eqnarray}
where $C_H>0$ and $\widehat{C}_H>0$ are two constants with
\[
C_H=\prod_{j=1}^dH_j(2H_j-1).
\]
Under assumption (\ref{intro-37}),
%
%
%e1.20 #&#
%
\begin{equation}
\label{intro-40} 0<\alpha\equiv2d-2\sum_{j=1}^dH_j<2
\wedge d.
\end{equation}

%
%th1.3 #&#
\begin{theorem}\label{intro-41} Assume (\ref{intro-37}).
For any $\theta>0$ and $x\in\R^d$,
%
%
%e1.21 #&#
%
\begin{eqnarray}
\label{intro-42} &&\lim_{t\to\infty}t^{-1}(\log
t)^{-{2}/{(4-\alpha)}}\log\E_x\exp \biggl\{ \theta \int
_0^t\frac{\partial^dW^H}{\partial x_1\cdots\partial
x_d}(B_s)\,ds
\biggr\}
\nonumber
\\[-8pt]
\\[-8pt]
\nonumber
&&\qquad=\frac{4-\alpha}{4} \biggl(\frac{\alpha}{2} \biggr)^{{\alpha}/ {(4-
\alpha)}} \bigl(2d C_H\theta^2\tilde{\kappa}(d,H)
\bigr)^{{2} /{(4-\alpha)}},\qquad \mbox{a.s.},
\nonumber
\end{eqnarray}
where
$\tilde{\kappa}(d, H)$ is the best constant of the inequality
[see (\ref{a-20}) in the \hyperref[app]{Appendix}]
\begin{eqnarray}
\int_{\R^d\times\R^d}f^2(x)f^2(y) \Biggl(\prod
_{j=1}^d\vert x_j-y_j
\vert^{2-2H_j} \Biggr)^{-1}\,dx\,dy\le C\|f\|_2^{4-\alpha}
\|\nabla f\| _2^\alpha, \nonumber\\
 \eqntext{f\in W^{1,2}\bigl(
\R^d\bigr).}
\end{eqnarray}
\end{theorem}

In the next theorem, we take $d=1$.
The Gaussian potential is a white noise formally
given as $V(x)=\dot{W}(x)$
where $W(x)$ $(x\in\R)$ is a two-side
Brownian motion. The relevant generalized Gaussian field is defined as
%
%
%e1.22 #&#
%
\begin{equation}
\label{intro-11} \langle V, \varphi\rangle=\int_{-\infty}^\infty
\varphi(x)W(dx),\qquad \varphi\in{\cal S}(\R).
\end{equation}
In this case
the covariance function $\gamma(\cdot)=\delta_0(\cdot)$
is the Dirac function
and the spectral measure
$\mu(d\lambda)=d\lambda$ is Lebesgue measure on $(-\infty, \infty)$.

%
%th1.4 #&#
\begin{theorem}\label{intro-46}
For any $\theta>0$ and $x\in\R$,
%
%
%e1.23 #&#
%
\begin{eqnarray}
\label{intro-47}&& \lim_{t\to\infty}t^{-1}(\log
t)^{-2/3}\log\E_x\exp \biggl\{\theta \int
_0^t \dot{W}(B_s)\,ds \biggr\}
\nonumber
\\[-8pt]
\\[-8pt]
\nonumber
&&\qquad =
\frac{1}{2} \biggl(\frac{3}{2} \biggr)^{2/3}
\theta^{4/3}, \qquad \mbox{a.s.}
\end{eqnarray}
\end{theorem}

We now comment on our main theorems.
It is interesting to see that (\ref{intro-35}) is consistent
with (\ref{intro-6}) when the latter is regarded as the case
$\alpha=0$, with easy and natural identifications $c(\gamma)=\gamma(0)$,
$\kappa(d,0)=1$, and the natural convention that $0^0=1$.

Given an integer
valued symmetric
simple
random walk $\{X_t;  t\ge0\}$ and
an independent family $\{\xi(x); x\in\Z\}$
of the i.i.d. standard normal random variables,
by Theorem~4.1, \cite{G-M}, or by Theorem~2.2, \cite{G-M-1},
\[
\lim_{t\to\infty}t^{-1}(\log t)^{-1/2}\log
\E_x\exp \biggl\{\theta \int_0^t
\xi(X_s)\,ds \biggr\} =\sqrt{2}\theta, \qquad \mbox{a.s.}
\]
Comparing this to Theorem~\ref{intro-46},
we witness a highly unusual difference between continuous and discrete
settings.

The almost sure limits stated in our theorems are largely determined
by the scaling or asymptotic scaling exponent $\alpha$
of the covariance function $\gamma(x)$
at $x=0$. The restriction $\alpha<2$
in our theorems is essential. In connection
to Theorem~\ref{intro-46}, notice that a Dirac
function on $\R^d$ satisfies $\delta_0(cx)=\vert c\vert^d\delta_0(x)$.
In particular, $\alpha=d$ as $\gamma(x)=\delta_0(x)$. To comply with
the restriction $\alpha<2$,
the space dimension $d$ has to be 1 in Theorem~\ref{intro-46}.

A challenge beyond the scope of this paper
is the quenched long-term asymptotics for the time dependent
potential $V(t,x)$ in connection to Theorems~\ref{intro-41}
and~\ref{intro-46}.
Associated to Theorem~\ref{intro-41} is the case
when
\[
V(t,x)=\frac{\partial^{d+1} W^H}{\partial t\,\partial x_1\cdots\partial
x_d}(t,x),
\qquad (t,x)\in\R^+\times\R^d,
\]
where $W^H(t, x)$
is a time--space
fractional Brownian sheet with some restriction on its Hurst parameter
$H=(H_0, H_1,\ldots, H_d)$.
An interested reader is referred to the paper by
Hu, Nualart and Song \cite{HNS} for the Feynman--Kac representation of the
solution in
this system; and to the recent work \cite{CHSX} by Chen, Hu, Song and Xing
for the annealed asymptotics
in this and other time--space settings.

Theorem~\ref{intro-46} corresponds to the famous
Karda--Parisi--Zhang (KPZ) model which starts from
a nonlinear stochastic partial differential equation
and is transformed into the parabolic Anderson equation
with the potential
\[
V(t, x)=\frac{\partial^2 W}{\partial t\,\partial x}(t,x),\qquad (t,x)\in\R^+\times\R
\]
by some renormalization treatment together with the Hopf--Cole
transform.
We cite the references \cite{karda-1} and \cite{karda-2} for
the physical background of the problem,
and \cite{ACQ,BQS,H} for the mathematical set-up and recent progress on the KPZ equation.

In addition, it is worth mentioning a recent
work \cite{CJKS} by Conus et al.
in which
they consider a possibly nonlinear heat equation
\[
{\partial}_t u=\tfrac{1}{2}\Delta u+V(t,x)\sigma(u).
\]
Here $V(t,x)$ is a time--space generalized Gaussian field with the covariance
function
\[
\delta_0(s-t)\gamma(x-y),\qquad (s,x), (t,y)\in\R^+\times\R^d.
\]
When the space covariance function $\gamma(x)$ satisfies (\ref{intro-15})
with $0<\alpha<2\wedge d$, a~quenched space-asymptotic law
(Theorem~2.6, \cite{CJKS}) states that
\[
C_1\le\limsup_{\vert x\vert\to\infty}\bigl(\log\vert x\vert
\bigr)^{-{2}/{
(4-\alpha)}} \log u(t,x)\le C_2,\qquad \mbox{a.s.}
\]
for any fixed $t>0$. The exponent $2/(4-\alpha)$ seems to suggest
a deep link to
(\ref{intro-35}). In general, going from the time-independent
potential to the time-dependent potential is a big step. We specially
mention the work \cite{VZ} by Viens and Zhang for their effort beyond
the sub-additivity treatment. It is our hope that
some ideas developed in the current paper may play a role in the future
investigation of this direction.

We now comment on the approaches adopted in this paper
and their relations to earlier works.
As usual, the proof consists of two major steps: a semi-group method
to associate the quenched exponential moment
in (\ref{intro-2}) to the principal eigenvalue of random linear operator
$2^{-1}\Delta+V$
with the zero boundary on $(-t,t)^d$
and asymptotic estimation of the principal eigenvalue for
which a nice idea developed in \cite{G-K} and \cite{G-K-M}
is adopted; see (\ref{bound-19}) to control
the principal eigenvalue
over the large domain $(-t,t)^d$
by the extreme among the principal eigenvalues
over the sub-domains. On the other hand, what sets this paper
apart is the singularity of our models. The following are some
of the novelties appearing in this paper.

\begin{longlist}[(1)]
\item[(1)] Algorithm development. The algorithms existing in the literature often
depend on the asymptotics of the generating function of $V(0)$.
Unfortunately, this strategy does
not apply here as $V(0)$ is not even defined in our models.
Indeed, the appearance of
$\|\nabla g\|_2$ in the constants of our main theorems
is a testimony of the dynamics different from the classic settings
represented by (\ref{intro-6}).
Our approach involves a rescaling strategy that highlights the
role of the diffusion part of the principal eigenvalue.
Some of the ideas adopted in this paper have been
used in the recent work \cite{C} in the setting of renormalized
Poissonian potential. However, there are substantial differences
between these two settings that demand some new adaptations.
The renormalized Poissonian potential is defined pointwise and essentially
total variational in the sense that it can be decomposed as
the difference of positive and negative parts
under suitable truncation,
while it is classic knowledge
that the potentials driven by white noise or fractional white noise
are not total variational.

\item[(2)] Entropy estimate. The entropy method has become an effective tool
in dealing with the tail, continuity,
integrability or finiteness for the random
quantities given as supremum.
In the case when $V(x)$ is defined pointwise, the concern is
the supremum $\sup_{x\in D}V(x)$ over a compact $D\subset\R^d$,
and the problem is
to count the $\varepsilon$-balls that cover $D$.
Not surprisingly, the entropy
number is bounded by a polynomial of $\varepsilon^{-1}$
if the distance is Euclidean or nearly Euclidean.
On the other hand, the entropy method in the context of generalized potential
is for the supremum $\sup_{g\in{\cal G}_d(D)}\langle V, g^2\rangle$
over (a dense set of)
the unit sphere of the Sobolev space over the domain $D$; see
Proposition~\ref{bound-3}. Counting the covering $\varepsilon$-balls
in a functional space is much harder and the result is less predictable
due to complexity in geometric structure.

\item[(3)] Lower bound by Slepian lemma. In the classic setting,
the lower bound for~(\ref{intro-6}) can be established
by decomposing $V(x)$ into two homogeneous Gaussian fields
such that the first field has finite correlation radius
and the second is negligible.
Under the assumptions of Theorems~\ref{intro-33} and~\ref{intro-41},
such decomposition is not available. Our treatment is based on
a famous comparison lemma by
Slepian~\cite{Slepian} and is formulated in Lemma~\ref{l-14} below.
\end{longlist}
%s2 #&#
\section{Gaussian supremum}\label{bound}

Let $D\subset\R^d$ be a fixed bounded open domain. We use the notation
${\cal S}(D)$ for the space of the infinitely smooth functions
on $D$ that vanish at the boundary of $D$.
For convenience,
we always view ${\cal S}(D)$ as a subspace of ${\cal S}(\R^d)$ by
defining $g(x)=0$ outside $D$ for each $g\in{\cal S}(D)$.
Given $g\in{\cal S}(D)$, for example, we may alternate between the notation
\[
\int_D\bigl\vert\nabla g(x)\bigr\vert^2\,dx \quad\mbox{and}\quad
\int_{\R^d}\bigl\vert\nabla g(x)\bigr\vert^2\,dx
\]
according to convenience. The notation $\|\nabla g\|_2$ is used for
both spaces ${\cal S}(D)$ and ${\cal S}(\R^d)$.
Set
%
%
%e2.1 #&#
%
\begin{eqnarray}
\label{bound-0} {\cal F}_d(D)& =& \bigl\{g\in{\cal S}(D); \|g
\|_2^2 =1 \bigr\},
\\
\label{bound-1} {\cal G}_d(D) &=& \bigl\{g\in{\cal S}(D); \|g
\|_2^2+\tfrac{1}{2}\|\nabla g\|_2^2
=1 \bigr\}.
\end{eqnarray}
Our approach largely relies on the estimate of
the supremum
%
%
%e2.3 #&#
%
\begin{equation}
\label{bound-2} \sup_{g\in{\cal F}_d(D)} \biggl\{\bigl\langle V,
g^2\bigr\rangle -\frac{1}{2}\int_{\R^d}\bigl\vert
\nabla g(x)\bigr\vert^2 \,dx \biggr\}.
\end{equation}
Notice that
for each $g\in{\cal S}(D)$,
$g^2\in{\cal S}(D)$. Consequently, the random variable $\langle V,
g^2\rangle$
is well defined and normal.
On the other hand, it is not obvious whether or not
the supremum is finite. When it is finite, the variation in (\ref{bound-2})
is the principal eigenvalue of the linear operator $(1/2)\Delta+V$
with the zero boundary condition over $D$.
The main goal of this section is
to show that the supremum in (\ref{bound-2}) is finite
when $D$ is bounded, and to establish a sharp almost-sure
asymptotic bound as $D$ expands to $\R^d$ in a suitable way.
The treatment is entropy estimation.

%s2.1 #&#
\subsection{Entropy bounds}

Consider a pseudometric space $(E, \rho)$ with the pseudometric
$\rho(\cdot, \cdot)$. For any $\varepsilon>0$, let $N(E,\rho,
\varepsilon)$
be the minimal number of the open balls of the diameter no greater than
$\varepsilon$,
which are necessary for covering $E$. In this section
we take $E={\cal G}_d(D)$ and
\[
\rho(f, g)= \bigl\{\E \bigl[\bigl\langle V, f^2\bigr\rangle-\bigl
\langle V, g^2\bigr\rangle \bigr]^2 \bigr\}^{1/2},
\qquad f,g\in{\cal G}_d(D).
\]
We have that
%
%
%e2.4 #&#
%
\begin{eqnarray}
\label{bound-2'} \qquad\rho(f, g)= \biggl\{\int_{\R^d\times\R^d}
\gamma(x-y) \bigl(f^2(x)-g^2(x) \bigr)
\bigl(f^2(y)-g^2(y) \bigr)\,dx\,dy \biggr\}^{1/2},
\nonumber
\\[-8pt]
\\[-8pt]
\eqntext{ f,g
\in{\cal G}_d(D).}
\end{eqnarray}
Here we specially mention that
$\gamma(x)=\delta_0(x)$ in the context of Theorem~\ref{intro-46}.

%
%pr2.1 #&#
\begin{proposition}\label{bound-3} Under the assumptions of
Theorems~\ref{intro-33}, \ref{intro-41} or~\ref{intro-46},
%
%
%e2.5 #&#
%
\begin{equation}
\label{bound-4} \lim_{\varepsilon\to0^+}\varepsilon^\beta \log N
\bigl({\cal G}_d(D), \rho, \varepsilon \bigr)=0
\end{equation}
whenever
%
%
%e2.6 #&#
%
\begin{equation}
\label{bound-4'} \beta> 1\vee\frac{2d}{ d+2}.
\end{equation}

Noticing that the right-hand side of (\ref{bound-4'}) is less than 2,
%
%
%e2.7 #&#
%
\begin{equation}
\label{bound-5} \int_0^1\sqrt{\log N \bigl({
\cal G}_d(D), \rho, \varepsilon \bigr)}\,d\varepsilon <\infty.
\end{equation}
\end{proposition}

\begin{pf}
Let $l(x)\in{\cal S}(\R^d)$ (mollifier) be a symmetric probability
density function supported on $\{\vert x\vert\le1\}$
and introduce the function
$l_\varepsilon(x)$ ($\varepsilon$-mollifier) as
%
%
%e2.8 #&#
%
\begin{equation}
\label{app-2} l_\varepsilon(x)=\varepsilon^{-d}l\bigl(
\varepsilon^{-1}x\bigr),\qquad x\in\R^d, \varepsilon>0.
\end{equation}
In addition, we assume that
${\cal F}(l)(\cdot)\ge0$.
Define the operator ${\cal S}_\varepsilon$ on ${\cal S}(\R^d)$
as
%
%
%e2.9 #&#
%
\begin{equation}
\label{bound-6} {\cal S}_\varepsilon g(x)= \biggl\{\int_{\R^d}g^2(x-y)l_\varepsilon
(y)\,dy \biggr\}^{1/2}, \qquad x\in\R^d.
\end{equation}
By Fourier transform,
\begin{eqnarray}
\E\bigl[\bigl\langle V, g^2\bigr\rangle-\bigl\langle V, {\cal
S}_\varepsilon (g)^2\bigr\rangle \bigr]^2 =
\frac{1}{(2\pi)^d}\int_{\R^d} \bigl\vert1-{\cal F}(l) (\varepsilon
\lambda) \bigr\vert^2 \bigl\vert{\cal F}\bigl(g^2\bigr) (\lambda)
\bigr\vert^2\mu(d\lambda),\nonumber\\
\eqntext{ g\in{\cal G}_d(D).}
\end{eqnarray}
Notice that $ \vert1-{\cal F}(l)(\varepsilon\lambda) \vert\le2$.
By the mean-value theorem there is $C_\delta>0$ such that
\[
\bigl\vert1-{\cal F}(l) (\varepsilon\lambda) \bigr\vert \le2^{1-\delta}\bigl \vert1-{\cal
F}(l) (\varepsilon\lambda) \bigr\vert ^\delta \le
C_\delta\vert\varepsilon
\lambda\vert^\delta,\qquad \lambda \in\R^d, \varepsilon>0,
\]
where $0<\delta<1$ is chosen by (\ref{intro-50}), in connection
to Lemma~\ref{intro-25'} in the \hyperref[app]{Appendix}.

Thus, there is a constant $C>0$ independent of $\varepsilon$ and $g$,
such that
\[
\rho (g, {\cal S}_\varepsilon g ) \le C\varepsilon^\delta \biggl\{\int
_{\R^d}\vert\lambda\vert ^{2\delta}\bigl \vert{\cal F}
\bigl(g^2\bigr) (\lambda) \bigr\vert^2\mu(d\lambda) \biggr
\}^{1/2},\qquad g\in{\cal G}_d(D), \varepsilon>0.
\]
Notice that
\[
\bigl\vert{\cal F}\bigl(g^2\bigr) (\lambda) \bigr\vert\le{\cal F}
\bigl(g^2\bigr) (0) =\|g\|_2^2\le1,\qquad g\in{\cal
G}_d(D).
\]
In addition, for any $\lambda\in\R^d\setminus\{0\}$,
\[
{\cal F}\bigl(g^2\bigr) (\lambda)=\frac{i}{ d} \int
_{\R^d} \biggl(\frac{\lambda}{\vert\lambda\vert^2}\cdot\nabla g^2(x)
\biggr) e^{i\lambda\cdot x}\,dx.
\]
Hence,
\begin{eqnarray*}\bigl\vert{\cal F}\bigl(g^2\bigr) (\lambda)
\bigr\vert&\le& \frac{1}{ d}\vert\lambda\vert^{-1}\int
_{\R^d}\bigl\vert\nabla g^2(x)\bigr\vert \,dx =
\frac{2}{ d}\vert\lambda\vert^{-1} \int_{\R^d}
\bigl\vert g(x)\bigr\vert\bigl\vert\nabla g(x)\bigr\vert \,dx
\\
&\le&\frac{2}{ d}\vert\lambda\vert^{-1}\|g\|_2\|
\nabla g\|_2\le \frac{2}{ d} \vert\lambda\vert^{-1}.
\end{eqnarray*}
Consequently,
\[
\bigl\vert{\cal F}\bigl(g^2\bigr) (\lambda) \bigr\vert^2\le
\biggl(1+\frac{2}{
d} \biggr) \biggl(1\wedge\frac{1}{\vert\lambda\vert^2} \biggr),\qquad g\in
{\cal G}_d(D).
\]
By (\ref{intro-50}), this leads to
\[
\sup_{g\in{\cal G}_d(D)}\int_{\R^d}\vert\lambda
\vert^{2\delta} \bigl\vert{\cal F}\bigl(g^2\bigr) (\lambda)
\bigr\vert^2\mu(d\lambda)<\infty.
\]

Summarizing our argument, there is a constant $C>0$ such that
%
%
%e2.10 #&#
%
\begin{equation}
\label{rho} \sup_{g\in{\cal G}_d(D)}\rho (g, {\cal S}_\varepsilon g )
\le C\varepsilon^\delta,\qquad \varepsilon>0.
\end{equation}
Write $\phi(\varepsilon)=\varepsilon^{\delta^{-1}}$. We have that
\[
\sup_{g\in{\cal G}_d(D)}\rho (g, {\cal S}_{\phi(\varepsilon)} g ) \le
C\varepsilon,\qquad \varepsilon>0.
\]

To prove (\ref{bound-4}), therefore, all we need is to show that
for any $\beta$ satisfying (\ref{bound-4'}),
%
%
%e2.11 #&#
%
\begin{equation}
\label{bound-7} \lim_{\varepsilon\to0^+}\varepsilon^\beta \log N
\bigl({\cal G}_d(D), \rho_\varepsilon, \varepsilon \bigr)=0,
\end{equation}
where the pseudometric $\rho_\varepsilon$ is defined as
$\rho_\varepsilon(f, g)
=\rho ({\cal S}_{\phi(\varepsilon)}f,  {\cal S}_{\phi
(\varepsilon
)}g )$
($f, g\in{\cal G}_d(D)$).
By (\ref{bound-2'})
\begin{eqnarray*} \rho_\varepsilon(f,g)&\le &\biggl(\int _{\R^d}\bigl \vert ({\cal
S}_{\phi(\varepsilon)}f)^2(x)-({ \cal S}_{\phi(\varepsilon)}g)^2(x)
\bigr\vert \,dx \biggr)^{1/2}
\\
&&{}\times \biggl(\sup_{x\in D'} \biggl\vert\int_{\R^d} \gamma(x-y)
\bigl\{ ({\cal S}_{\phi(\varepsilon)}f)^2(y)-({\cal
S}_{\phi(\varepsilon )}g)^2(y) \bigr\} \,dy \biggr\vert \biggr)^{1/2},
\end{eqnarray*}
where $D'$ is the 1-neighborhood of $D$. Take
\begin{eqnarray}
A_\varepsilon(f) (x)=({\cal S}_{\phi(\varepsilon)}f)^2(x)\quad \mbox
{and}\quad B_\varepsilon(f) (x)=\int_{\R^d}\gamma(x-y) ({\cal
S}_{\phi(\varepsilon)}f)^2(y)\,dy,\nonumber\\
 \eqntext{x\in D'}
\end{eqnarray}
in Lemma~\ref{ball-3} of the \hyperref[app]{Appendix}. All we need is to exam that
there are $p>1$
satisfying
%
%
%e2.12 #&#
%
\begin{equation}
\label{bound-11'} \beta>\frac{2p}{2p-1}>1\vee\frac{2d}{ d+2}
\end{equation}
and $C>0$, $m>0$ independent of $\varepsilon>0$ such that
%
%
%e2.13 #&#
%
\begin{eqnarray}
\label{bound-13} \bigl\vert({\cal S}_{\phi(\varepsilon)}f)^2(x)- ({\cal
S}_{\phi(\varepsilon)}f)^2(y) \bigr\vert&\le& C\varepsilon^{-m}\vert
x-y \vert,
\\
\label{bound-12} \biggl\vert \int_{\R^d} \bigl\{\gamma(x-z)-
\gamma(y-z) \bigr\} ({\cal S}_{\phi(\varepsilon)}f)^2(z)\,dz
\biggr\vert&\le& C \varepsilon ^{-m}\vert x-y\vert,
\\
\label{bound-15'} \int_{\R^d}\bigl\vert({\cal S}_{\phi(\varepsilon)}f)
(z)\bigr\vert^{2p}\,dz&\le& C
\end{eqnarray}
and
%
%
%e2.16 #&#
%
\begin{equation}
\label{bound-14} \biggl\vert \int_{\R^d}\gamma(x-z) ({\cal
S}_{\phi(\varepsilon)}f)^2(z)\,dz \biggr\vert\le C
\end{equation}
for all $x, y\in D'$ and $f\in{\cal G}_d(D)$.

Indeed, by the mean value theorem
\begin{eqnarray*} \bigl \vert({\cal S}_{\phi(\varepsilon)}f)^2(x)- ({\cal
S}_{\phi(\varepsilon)}f)^2(y)\bigr \vert &\le&\int_{\R^d} \bigl\vert
l_{\phi(\varepsilon)}(x+z)-l_{\phi(\varepsilon )}(y+z)\bigr \vert
f^2(z)\,dz
\\
&\le &C\phi(\varepsilon)^{-(d+1)}\vert x-y\vert\int_{\R^d}
f^2(z)\,dz\\
&\le& C\phi(\varepsilon)^{-(d+1)}\vert x-y\vert.
\end{eqnarray*}
Thus (\ref{bound-13}) holds with $m=(d+1)\delta^{-1}$.
For the same $m$,
(\ref{bound-12}) follows from~(\ref{bound-13}), the relation
\begin{eqnarray*}&&\int_{\R^d} \bigl\{\gamma(x-z)- \gamma(y-z) \bigr\} ({\cal
S}_{\phi(\varepsilon)}f)^2(z)\,dz
\\
&&\qquad=\int_{\R^d}\gamma(z) \bigl\{({\cal
S}_{\phi(\varepsilon)}f)^2(z-x) -({\cal S}_{\phi(\varepsilon)}f)^2(z-y)
\bigr\}\,dz,
\end{eqnarray*}
and the fact that
\[
\int_{\widetilde{D}}\bigl\vert\gamma(z)\bigr\vert \,dz<\infty
\]
for $\widetilde{D}=\{z_1+z_2\in\R^d;  z_1, z_2\in D'\}$.

We now come to (\ref{bound-15'}). First, for any $p>1$ and
by Jensen's inequality,
\[
\int_{\R^d}\bigl\vert({\cal S}_{\phi(\varepsilon)}f) (z)\bigr \vert^{2p}\,dz
\le\int_{\R^d}\bigl\vert f(z) \bigr\vert^{2p}\,dz.
\]
We claim that there is a $p>1$ satisfying (\ref{bound-11'}) and
$p(d-2)<d$. Indeed, this is obvious when $d\le2$ as we can make $p$
sufficiently large. When $d\ge3$, our assertion is secured
by the facts that the quantity $2p(2p-1)^{-1}$ is
strictly decreasing in~$p$, and that
the supremum of $p$
under the constraint $p(d-2)<d$ is $b\equiv d(d-2)^{-1}$ which
solves the equation
\[
\frac{2b}{2b-1}=\frac{2d}{ d+2}.
\]
By Gagliardo--Nirenberg inequality
(see, e.g., page 303, \cite{Chen}), for which the restriction
$p(d-2)<d$ is critically needed,
\[
\int_{\R^d}\bigl \vert f(x) \bigr\vert^{2p}\,dx \le C\|f
\|_2^{d(p-1)}\|\nabla f\|_2^{2p-d(p-1)}\le C.
\]
Thus, we have proved (\ref{bound-15'}).

It remains to establish (\ref{bound-14}).
In the context of Theorem~\ref{intro-41},
by (\ref{a-18}),
\[
\biggl\vert \int_{\R^d}\gamma(x-z) ({\cal S}_{\phi(\varepsilon)}f)^2(z)\,dz
\biggr\vert\le C\|{\cal S}_{\phi(\varepsilon)}f\|_2^{4-\alpha} \|\nabla{
\cal S}_{\phi(\varepsilon)}f\|_2^{\alpha}.
\]
By Jensen inequality, $\|{\cal S}_{\phi(\varepsilon)}f\|_2\le\|f\|_2\le1$.
From (\ref{bound-6})
\begin{eqnarray*}&&\bigl\vert\nabla{\cal S}_{\phi(\varepsilon)}f(x)\bigr\vert
\\
&&\qquad=
\biggl(\int_{\R^d}l_{\phi(\varepsilon)}(y) f^2(x-y)\,dy
\biggr)^{-1/2} \biggl\vert\int_{\R^d}l_{\phi(\varepsilon)}(y)f(x-y)
\nabla f(x-y)\,dy \biggr\vert
\\
&&\qquad\le \biggl\{\int_{\R^d}l_{\phi(\varepsilon)}(y)\bigl\vert\nabla
f(x-y)\bigr\vert^2\,dy \biggr\}^{1/2},
\end{eqnarray*}
where the inequality follows from Cauchy--Schwarz
inequality. Hence, by Fubini's theorem and translation invariance,
%
%
%e2.17 #&#
%
\begin{equation}
\label{bound-14'} \|\nabla{\cal S}_{\phi(\varepsilon)}f
\|_2^2\le\int_{\R^d}l_{\phi
(\varepsilon)}(y)
\biggl[\int_{\R^d}\bigl\vert\nabla f(x-y)\bigr\vert^2\,dx
\biggr]\,dy =\|\nabla f\|_2^2.
\end{equation}
The right-hand side is bounded by 1.
Thus (\ref{bound-14}) holds.

In the context of Theorem~\ref{intro-33},
(\ref{bound-14}) follows from the
bound
$\vert\gamma(z)\vert\le C(1+\vert z\vert^{-\alpha})$ and a similar estimate
[with (\ref{a-18}) being replaced by (\ref{a-1})].

In the context of Theorem~\ref{intro-46},
\[
\int_{\R^d}\gamma(x-z) ({\cal S}_{\phi(\varepsilon)}f)^2(z)\,dz=({
\cal S}_{\phi(\varepsilon)}f)^2(x) \le\sup_{y\in\R}f^2(y).
\]
Hence, (\ref{bound-14}) follows from the estimate
\begin{eqnarray*}
f^2(y)&\le&2\int_{-\infty}^\infty\bigl\vert f(u)
f'(u)\bigr\vert \,du\\
&\le& 2 \biggl\{\int_{-\infty}^\infty
f^2(u)\,du \biggr\}^{1/2} \biggl\{\int_{-\infty}^\infty
\bigl\vert f'(u)\bigr\vert^2\,du \biggr\}^{1/2}\le 2,\qquad y\in
\R.
\end{eqnarray*}
\upqed\end{pf}

%s2.2 #&#
\subsection{Consequences of the entropy bounds}

According to the classic theory on sample path regularity (see, e.g.,
Appendix D, \cite{Chen}),
under (\ref{bound-5})
the supremum in~(\ref{bound-2}) is finite,
integrable and $\{\langle V, g^2\rangle;  g\in{\cal
G}_d(D)\}
$ has
continuous sample paths with respect to
the pseudometric induced by its covariance. By the linearity
of $V$ and a standard extension argument, such sample continuity
is extended to ${\cal S}(\R^d)$.

Given a generalized function $\xi$ on $\R^d$, that is, a linear functional
on ${\cal S}_d(\R^d)$, set
%
%
%e2.18 #&#
%
\begin{equation}
\label{bound-16} \lambda_\xi(D)=\sup_{g\in{\cal F}_d(D)} \biggl
\{\bigl\langle\xi, g^2\bigr\rangle- \frac{1}{2}\int
_{D}\bigl\vert\nabla g(x)\bigr\vert^2 \,dx \biggr\}.
\end{equation}
For any $\varepsilon>0$, let $D_\varepsilon$ be the $\varepsilon$-neighborhood
of $D$. By the obvious monotonicity of $\lambda_\xi(D)$ in $D$, the limit
%
%
%e2.19 #&#
%
\begin{equation}
\label{bound-16'} \lambda_\xi^+(D)\equiv\lim
_{\varepsilon\to0^+}\lambda_\xi (D_\varepsilon)
\end{equation}
always exists at least as extended number. It is not clear to us
whether or when $\lambda_\xi^+(D)=\lambda_\xi(D)$.

Let the $\varepsilon$-mollifier $l_\varepsilon(\cdot)$ be given in (\ref{app-2})
and define the pointwise random field $V_\varepsilon(\cdot)$ as
%
%
%e2.20 #&#
%
\begin{equation}
\label{app-3} V_\varepsilon(x)=\bigl\langle V, l_\varepsilon(\cdot-x)
\bigr\rangle,\qquad x\in \R^d.
\end{equation}

%
%le2.2 #&#
\begin{lemma}\label{bound-14''} Under the assumptions
of Theorems~\ref{intro-33},
\ref{intro-41} or~\ref{intro-46}
%
%
%e2.21 #&#
%
\begin{equation}
\label{lemma-1} \lim_{\varepsilon\to0^+}\E\sup_{g\in{\cal G}_d((-\varepsilon, \varepsilon)^d)}
\bigl\langle V, g^2\bigr\rangle=0
\end{equation}
and
%
%
%e2.22 #&#
%
\begin{equation}
\label{lemma-2} \lambda_{\theta V}(D)\le\liminf_{\varepsilon\to0^+}
\lambda_{\theta
V_\varepsilon}(D)\le \limsup_{\varepsilon\to0^+}
\lambda_{\theta V_\varepsilon}(D)\le\lambda _{\theta V}^+(D), \qquad \mbox{a.s.}
\end{equation}
for any $\theta>0$ and bound domain $D\subset\R^d$.
\end{lemma}

\begin{pf} In our view, ${\cal G}_d((-\varepsilon, \varepsilon)^d)$ is a subset
of ${\cal G}_d ((-1, 1)^d )$ as $\varepsilon<1$.
By the continuity of the Gaussian field
$\{\langle V, g^2\rangle;  g\in{\cal G}_d ((-1,
1)^d )\}$
with respect to its covariance function
established by Proposition~\ref{bound-3},
\[
\lim_{\delta\to0^+}\E\sup \bigl\{\bigl\langle V, g^2\bigr
\rangle; g\in{\cal G}_d \bigl((-1, 1)^d \bigr) \mbox{ and }
\E\bigl\langle V, g^2\bigr\rangle^2\le\delta \bigr\}=0.
\]
To establish (\ref{lemma-1}), it suffices to examine that
%
%
%e2.23 #&#
%
\begin{equation}
\label{bound-15} \lim_{\varepsilon\to0^+}\sup_{g\in{\cal G}_d((-\varepsilon, \varepsilon)^d)} \E
\bigl\langle V, g^2\bigr\rangle^2=0.
\end{equation}

Indeed, in the case of Theorem~\ref{intro-33},
\begin{eqnarray*}\E\bigl\langle V, g^2\bigr
\rangle^2&=&\int_{\R^d}\gamma(x-y)
g^2(x)g^2(y)\,dx\,dy
\\
&\le& C\int_{\R^d}\frac{g^2(x)g^2(y)}{\vert x-y\vert^\alpha}\,dx\,dy \le C
\varepsilon^{\alpha'-\alpha} \int_{\R^d}\frac{g^2(x)g^2(y)}{\vert
x-y\vert^{\alpha'}}\,dx\,dy,
\end{eqnarray*}
where the constant $C>0$ is different in each step but independent of
$g$. The constant $\alpha'$ is chosen by the principle that
$\alpha<\alpha'< 2\wedge d$.
Consequently,
\[
\int_{\R^d}\frac{g^2(x)g^2(y)}{\vert x-y\vert^{\alpha'}}\,dx\,dy \le C_{\alpha'}\|g
\|_2^{2-\alpha'}\vert\nabla g\|_2^{\alpha'}\le
C_{\alpha'}, \qquad g\in{\cal G}_d\bigl((-\varepsilon,
\varepsilon)^d\bigr),
\]
where $C_{\alpha'}$ is given in (\ref{a-2})
with $\alpha$ being replaced by $\alpha'$. Hence, we have
(\ref{bound-15}).

This argument applies also to the settings of Theorems~\ref{intro-41} and~\ref{intro-46}.
For Theorem~\ref{intro-41}, we use (\ref{a-18})
instead of (\ref{a-1}) and pick $2H_j-1 <\alpha_j<1$ ($j=1,\ldots, d$)
with $\alpha_1+\cdots+\alpha_d<2$.

As for Theorem~\ref{intro-46}, we first apply
in (A.2), \cite{BCR} [with $p=d=1$, $\sigma=1/2$
and $f(x)=g^4(x)$] that gives
\[
\int_{-\infty}^\infty\frac{g^4(x)}{\vert x\vert^{1/2}}\,dx \le C\|g
\|_{8}^4,\qquad g\in{\cal G}_d(\R),
\]
where $C>0$ is independent of $g$.
The right-hand side is uniformly bounded over $g\in{\cal G}_d(\R)$
according to the Gagliardo--Nirenberg
inequality (see, e.g., (C.1), page~303, \cite{Chen})
\[
\|g\|_{8}\le C\bigl\|g'\bigr\|_2^{{3}/{8}}\|g\|^{{5} /8} \le C,\qquad g\in{\cal G}_d(\R).
\]

We now come to (\ref{lemma-2}). Let $g\in{\cal F}_d(D)$ be fixed but
arbitrary.
\[
\lambda_{\theta V_\varepsilon}(D)\ge\theta\int_{\R^d}
V_\varepsilon (x)g^2(x)\,dx-\frac{1}{2} \int
_{\R^d}\bigl\vert\nabla g(x)\bigr\vert^2 \,dx.
\]
By linearity,
%
%
%e2.24 #&#
%
\begin{equation}
\label{lemma-3} \int_{\R^d} V_\varepsilon(x)g^2(x)\,dx=
\int_{\R^d} \bigl\langle V, l_\varepsilon(\cdot-x)\bigr
\rangle g^2(x)\,dx =\bigl\langle V, ({\cal S}_\varepsilon
g)^2\bigr\rangle.
\end{equation}
In addition, by (\ref{rho}) and a proper normalization
one can see
that ${\cal S}_\varepsilon g$ converges to $g$ under the covariance
pseudomatric
$\rho$ given in (\ref{bound-2'}). By the sample path continuity of
the functional $\langle V, g^2\rangle$ resulting from
Proposition~\ref{bound-3},
\[
\lim_{\varepsilon\to0^+}\bigl\langle V, ({\cal S}_\varepsilon
g)^2\bigr\rangle =\bigl\langle V, g^2\bigr\rangle,\qquad \mbox{a.s.}
\]
Hence,
\[
\liminf_{\varepsilon\to0^+}\lambda_{\theta V_\varepsilon}(D)\ge \theta\bigl\langle
V, g^2\bigr\rangle-\frac{1}{2} \int_{\R^d}
\bigl\vert\nabla g(x)\bigr\vert^2 \,dx,\qquad \mbox{a.s.}
\]
Taking supremum over $g$ on the right-hand side, we establish
the lower bound needed by (\ref{lemma-2}).

As for the upper bound, first notice that for any $g\in{\cal F}_d(D)$,
$f\equiv\break \|{\cal S}_\varepsilon g\|_2^{-1}{\cal S}_\varepsilon g\in
{\cal F}_d(D_\varepsilon)$.
By (\ref{lemma-3}) and linearity,
\begin{eqnarray*} \lambda_{\theta V_\varepsilon}(D)&\le&\sup
_{g\in{\cal F}_d(D)} \biggl\{\bigl\langle V, ({\cal S}_\varepsilon
g)^2\bigr\rangle-\frac{1}{2}\int_{\R^d}
\bigl\vert\nabla({\cal S}_\varepsilon g) (x)\bigr\vert^2 \,dx \biggr\}
\\
&\le& \Bigl(\sup_{g\in{\cal F}_d(D)}\|{\cal S}_\varepsilon g
\|_2^2 \Bigr) \sup_{f\in{\cal F}_d(D_\varepsilon)} \biggl\{\bigl
\langle V, f^2\bigr\rangle-\frac{1}{2}\int_{\R^d}
\bigl\vert\nabla f(x)\bigr\vert^2 \,dx \biggr\} \\
&\le&\lambda_{\theta V}(D_\varepsilon),
\end{eqnarray*}
where the last step follows from the fact
$\|{\cal S}_\varepsilon g\|_2\le\|g\|_2 =1$ [see (\ref{bound-14'})]
for any $g\in{\cal F}_d(D)$.

Letting $\varepsilon\to0^+$ leads to the upper bound needed by (\ref{lemma-2}).
\end{pf}

In the rest of the section, we demonstrate how Proposition~\ref{bound-3}
(or Lem\-ma~\ref{bound-14''}, more precisely)
is used to bound the principal eigenvalue given in (\ref{bound-2}).

The principal eigenvalue over a large domain can be essentially bounded
by the extreme value among the principal eigenvalues of the sub-domains,
according to a nice strategy developed by
G\"artner and K\"onig \cite{G-K}.
Let $r\ge2$. By Proposition~1 in \cite{G-K}, also by
Lemma~4.6 in \cite{G-K-M}, there is
a nonnegative and continuous
function $\Phi(x)$ on $\R^d$
whose support is contained in the 1-neighborhood of the grid $2r\Z^d$,
such that for any $R>r$ and any generalized function $\xi$,
%
%
%e2.25 #&#
%
\begin{equation}
\label{bound-17} \lambda_{\xi-\Phi^y}(Q_R)\le\max
_{z\in2r\Z^d\cap Q_R}\lambda _\xi(z+Q_{r+1}),\qquad y\in
Q_r,
\end{equation}
where $\Phi^y(x)=\Phi(x+y)$, and we use the notation
$Q_R=(-R, R)^d$ for any $R>0$.

In addition, $\Phi(x)$ is periodic with period $2r$,
\[
\Phi(x+2rz)=\Phi(x),\qquad x\in\R^d, z\in\Z^d,
\]
and there is a constant $K>0$ independent of $r$ such that
%
%
%e2.26 #&#
%
\begin{equation}
\label{bound-18} \frac{1}{(2r)^d}\int_{Q_{r}}\Phi(x)\,dx\le
\frac{K}{ r}.
\end{equation}
It should be pointed out that originally, (\ref{bound-17}) was established
for the ordinary function $\xi$. However, it can be extended to
the generalized function without any extra effort, due to the linearity
preserved by the form $\langle\xi,\varphi\rangle$
($\varphi\in{\cal S}(\R^d)$).

Write
\[
\eta(x)=\frac{1}{(2r)^d}\int_{Q_r}\Phi(x+y)\,dy=
\frac{1}{(2r)^d}\int_{Q_r}\Phi^y(x)\,dy,\qquad x\in
\R^d.
\]
By periodicity, $\eta\equiv\eta(x)$ is a constant with a bound
given in (\ref{bound-18}).
Hence,
%
%
%e2.27 #&#
%
\begin{eqnarray}
\label{bound-19} \lambda_{\xi}(Q_R)&\le&
\frac{K}{ r}+\lambda_{\xi-\eta}(Q_R) \le
\frac{K}{ r}+\frac{1}{(2r)^d}\int_{Q_r}
\lambda_{\xi-\Phi
^y}(Q_R)\,dy
\nonumber
\\[-8pt]
\\[-8pt]
\nonumber
&\le&\frac{K}{ r}+\max_{z\in2r\Z^d\cap Q_R}\lambda_\xi
(z+Q_{r+1}),
\end{eqnarray}
where the last inequality follows from (\ref{bound-17}), and
the second inequality follows from the following steps:
\begin{eqnarray*}\lambda_{\xi-\eta}(Q_R)&=&\sup
_{g\in{\cal F}_d(Q_R)} \biggl\{ \frac{1}{(2r)^d}\int_{Q_r}
\bigl\langle\xi-\Phi^y, g^2\bigr\rangle \,dy-
\frac{1}{2} \int_{Q_R}\bigl\vert\nabla g(x)\bigr\vert^{2}\,dx \biggr\}
\\
&=&\sup_{g\in{\cal F}_d(Q_R)} \biggl\{ \frac{1}{(2r)^d}\int
_{Q_r} \biggl[\bigl\langle\xi-\Phi^y,
g^2\bigr\rangle \,dy-\frac{1}{2} \int_{Q_R}
\bigl\vert\nabla g(x)\bigr\vert^{2}\,dx \biggr]\,dy \biggr\}
\\
&\le&\frac{1}{(2r)^d}\int_{Q_r}\sup_{g\in{\cal F}_d(Q_R)}
\biggl[\bigl\langle\xi-\Phi^y, g^2\bigr\rangle \,dy-
\frac{1}{2} \int_{Q_R}\bigl\vert\nabla g(x)\bigr\vert^{2}\,dx \biggr]\,dy
\\
&=&\frac{1}{(2r)^d}\int_{Q_r}\lambda_{\xi-\Phi^y}(Q_R)\,dy.
\end{eqnarray*}

In the next lemma, we not only show that the principal eigenvalue
in (\ref{bound-2}) is finite for any bounded domain $D$, but also
provide sharp asymptotic bounds
for the almost-sure increasing rate of the principal eigenvalue
as $D$ expands to $\R^d$ in a proper way.

%
%le2.3 #&#
\begin{lemma}\label{bound-20} Under the assumptions of Theorems~\ref{intro-33}
or~\ref{intro-41}, for any $\theta>0$,
%
%
%e2.28 #&#
%
\begin{equation}
\label{bound-21} \limsup_{t\to\infty}(\log t)^{-{2}/{(4-\alpha)}}
\lambda_{\theta V}(Q_t) \le\theta^{{4}/ {(4-
\alpha)}}h(d,\alpha),\qquad \mbox{a.s.},
\end{equation}
where
%
%
%e2.29 #&#
%
\begin{equation}
\label{bound-22} h(d,\alpha)=\cases{ %
\displaystyle\frac{4-\alpha}{4} \biggl(\frac{\alpha}{2} \biggr)^{{\alpha}/ { (4-
\alpha)}} \bigl(2\,dc(\gamma)\kappa(d,\alpha) \bigr)^{{2}
/{(4-\alpha)}},\vspace*{2pt}\cr
\qquad \mbox{in the setting of Theorem~\ref{intro-33},}
\vspace*{2pt}\cr
\displaystyle\frac{4-\alpha}{4} \biggl(\frac{\alpha}{2} \biggr)^{{\alpha}/ { (4-
\alpha)}} \bigl(2dC_H\tilde{\kappa}(d,H) \bigr)^{{2}/ {(4-
\alpha)}},\vspace*{2pt}\cr
\qquad \mbox{in the setting of Theorem~\ref{intro-41}.}}
\end{equation}

Under the assumption of Theorem~\ref{intro-46}, for any $\theta>0$,
%
%
%e2.30 #&#
%
\begin{equation}
\label{bound-23} \limsup_{t\to\infty}(\log t)^{-2/3}
\lambda_{\theta V} \bigl((-t, t) \bigr) \le\frac{1}{2} \biggl(
\frac{3}{2} \biggr)^{2/3}\theta^{4/3},\qquad\mbox{a.s.}
\end{equation}
\end{lemma}

\begin{pf} Let $u>0$ be fixed, and write
%
%
%e2.31 #&#
%
\begin{equation}
\label{u-7} a(t)=\cases{ %
\sqrt{u}(\log
t)^{{1}/ {(4-\alpha)}},\vspace*{2pt}\cr \qquad\mbox{in the setting of Theorems~\ref{intro-33}
or~\ref{intro-41},}
\vspace*{2pt}\cr
\sqrt{u}(\log t)^{1/3},\vspace*{2pt}\cr \qquad\mbox{in the setting of
Theorem~\ref{intro-46}.}}
\end{equation}
For each $g\in{\cal S}(\R^d)$, write
%
%
%e2.32 #&#
%
\begin{equation}
\label{u-8} g_t(x)=a(t)^{d/2}g \bigl(a(t)x \bigr),\qquad x\in
\R^d.
\end{equation}
By rescaling substitution $g\mapsto g_t$,
%
%
%e2.33 #&#
%
\begin{equation}
\label{u-7'} \qquad\lambda_{\theta V}(Q_t)=a(t)^2
\sup_{g\in{\cal F}_d(Q_{ta(t)})} \biggl\{ \theta a(t)^{-2} \bigl\langle V,
g_t^2\bigr\rangle -\frac{1}{2}\int
_{Q_{ta(t)}}\bigl\vert\nabla g(x)\bigr\vert^2\,dx \biggr\}.
\end{equation}
Let
$\{\langle V_t,\varphi\rangle;
\varphi\in{\cal S}(\R^d)\}$ be the generalized
Gaussian field defined as $\langle V_t, \varphi)\rangle
=\langle V, \tilde{\varphi}_t\rangle$,
where $\tilde{\varphi}(x)=a(t)^d\varphi (a(t)x )$ [notice that
this is different from the definition in (\ref{u-8})].
Then we have $\langle V, g_t^2\rangle=\langle V_t, g^2\rangle$.
Taking $\xi=\theta a(t)^{-2}V_t$
in (\ref{bound-19}), by (\ref{u-7'}) we have that
%
%
%e2.34 #&#
%
\begin{equation}
\label{u-9} \lambda_{\theta V}(Q_t)\le a(t)^2 \biggl
\{\frac{K}{ r} + \max_{z\in2r\Z^d\cap Q_{ta(t)}}X_z(t) \biggr\}
\end{equation}
for any $r\ge2$, where, by homogeneity of the Gaussian field
$\{\langle V, \varphi\rangle;  \varphi\in{\cal S}(\R
^d)\}$,
the stochastic
processes
\begin{eqnarray}
X_z(t)\equiv\sup_{g\in{\cal F}_d(z+Q_{r+1})} \biggl\{\theta
a(t)^{-2} \bigl\langle V, g_t^2\bigr\rangle-
\frac{1}{2}\int_{z+Q_{r+1}}\bigl\vert\nabla g(x)
\bigr\vert^2\,dx \biggr\},\nonumber\\
 \eqntext{z\in2r\Z^d\cap Q_{ta(t)}}
\end{eqnarray}
are identically distributed.
Thus
\[
\P \Bigl\{\max_{z\in2r\Z^d\cap Q_{ta(t)}}X_z(t)> 1 \Bigr\} \le\#
\bigl\{2r\Z^d\cap Q_{ta(t)}\bigr\} \P\bigl
\{X_0(t)> 1\bigr\}.
\]

By linearity, for any $g\in{\cal F}_d(Q_{r+1})$,
\begin{eqnarray*} &&\theta a(t)^{-2} \bigl\langle V,
g_t^2\bigr\rangle-\frac{1}{2}\int
_{z+Q_{r+1}}\bigl\vert\nabla g(x)\bigr\vert ^2\,dx
\\
&&\qquad\le\theta a(t)^{-2} \Bigl(\sup_{f\in{\cal G}_d(Q_{r+1})}\bigl\langle
V, f_t^2\bigr\rangle \Bigr) \biggl(1+\frac{1}{2} \|
\nabla g\|_2^2 \biggr)-\frac{1}{2}\|\nabla g
\|_2^2.
\end{eqnarray*}
Here we recall that
the class ${\cal G}_d(D)$ is defined in (\ref{bound-1}).
Taking supremum over $g$,
\[
X_0(t)\le\sup_{g\in{\cal F}_d(Q_{r+1})} \biggl\{ \theta
a(t)^{-2} \Bigl(\sup_{f\in{\cal G}_d(Q_{r+1})}\bigl\langle V,
f_t^2\bigr\rangle \Bigr) \biggl(1+\frac{1}{2} \|
\nabla g\|_2^2 \biggr)-\frac{1}{2}\|\nabla g
\|_2^2 \biggr\}.
\]
Consequently,
\[
\bigl\{X_0(t)\ge1\bigr\}\subset \Bigl\{\sup_{f\in{\cal G}_d(Q_{r+1})}
\bigl\langle V, f_t^2\bigr\rangle\ge\theta^{-1}
a(t)^2 \Bigr\}.
\]

Summarizing our argument,
%
%
%e2.35 #&#
%
\begin{eqnarray}
\label{u-10} &&\P \Bigl\{\max_{z\in2r\Z^d\cap Q_{ta(t)}}X_z(t)> 1 \Bigr
\}
\nonumber
\\[-8pt]
\\[-8pt]
\nonumber
&&\qquad\le\# \bigl\{2r\Z^d\cap Q_{ta(t)}\bigr\}\P \Bigl\{\sup
_{g\in{\cal G}_d(Q_{r+1})}\bigl\langle V, g_t^2\bigr\rangle
\ge \theta^{-1} a(t)^2 \Bigr\}.
\end{eqnarray}

Notice that for each $g\in{\cal G}_d(Q_{r+1})$,
$ (1+a(t)^2\|\nabla g\|_2^2 )^{-1/2}g_t(\cdot)
\in\break  {\cal G}_d(Q_{(r+1)a(t)^{-1}})$.
By linearity,
\begin{eqnarray*}
\E\sup_{g\in{\cal G}_d
(Q_{r+1})}\bigl\langle V, g_t^2
\bigr\rangle&\le& \bigl(1+a(t)^2 \bigr) \E\sup_{f\in{\cal G}_d(Q_{(r+1)a(t)^{-1}})}
\bigl\langle V, f^2\bigr\rangle\\
& =&o \bigl(a(t)^2 \bigr)\qquad (t
\to\infty),
\end{eqnarray*}
where the last step follows from (\ref{lemma-1}) in
Lemma~\ref{bound-14''}.

By the concentration inequality for Gaussian field (see, e.g., (5.152),
Theorem~5.4.3, page 219, \cite{MR}, in connection to Corollary~5.4.5,
page 224,
\cite{MR}),
%
%
%e2.36 #&#
%
\begin{eqnarray}
\label{u-10'} &&\P \Bigl\{\sup_{g\in{\cal G}_d(Q_{r+1})}\bigl\langle
V, g_t^2\bigr\rangle >\theta^{-1}
a(t)^2 \Bigr\}
\nonumber\\
&&\qquad=\P \Bigl\{\sup_{g\in{\cal G}_d(Q_{r+1})}\bigl\langle V, g_t^2
\bigr\rangle -\E\sup_{g\in{\cal G}_d(Q_{r+1})}\bigl\langle V, g_t^2
\bigr\rangle > \bigl(1+o(1) \bigr)\theta^{-1} a(t)^2 \Bigr\}
\\
&&\qquad\le\exp \biggl\{- \bigl(1+o(1) \bigr)\frac{a(t)^4}{2\theta^{2}\sigma
_t^2} \biggr\},
\nonumber
\end{eqnarray}
where
\[
\sigma_t^2=\sup_{g\in{\cal G}_d(Q_{r+1})}\Var \bigl(
\bigl\langle V, g_t^2\bigr\rangle \bigr).
\]

In the setting of Theorem~\ref{intro-33}, by (\ref{intro-15})
and other assumptions on $\gamma(x)$,
\begin{eqnarray*}\sigma_t^2&=&\sup
_{g\in{\cal G}_d(Q_{r+1})} \int_{\R^d\times\R^d}\gamma(x-y)g_t^2(x)g_t^2(y)\,dx\,dy
\\
&=&\sup_{g\in{\cal G}_d(Q_{r+1})}\int_{\R^d\times\R^d} \gamma
\bigl(a(t)^{-1}(x-y) \bigr)g^2(x)g^2(y)\,dx\,dy
\\
&\sim& c(\gamma)a(t)^\alpha\sup_{g\in{\cal G}_d(Q_{r+1})}\int
_{\R
^d\times\R^d} \frac{g^2(x)g^2(y)}{\vert x-y\vert^\alpha}\,dx\,dy \qquad(t\to\infty).
\end{eqnarray*}
Notice that
\begin{eqnarray*}
\sup_{g\in{\cal G}_d(Q_{r+1})}\int_{\R^d\times\R^d}
\frac{g^2(x)g^2(y)}{\vert x-y\vert^\alpha}\,dx\,dy &\le&\sigma^2(d,\alpha)\\
&=& \biggl(\frac{4-\alpha}{4}
\biggr)^{{(4-\alpha)}/ {2} }\biggl(\frac{\alpha}{2}
\biggr)^{\alpha/2}\kappa(d,\alpha),
\end{eqnarray*}
where $\sigma(d, \alpha)$
is the variation defined in (\ref{a-5})
and the last step follows from~(\ref{a-8}) of Lemma~\ref{a-6} in the \hyperref[app]{Appendix}.

In view of (\ref{u-7}),
%
%
%e2.37 #&#
%
\begin{eqnarray}
\label{u-11} &&\P \Bigl\{\sup_{g\in{\cal G}_d(Q_{r+1})}\bigl\langle V,
g_t^2\bigr\rangle >\theta^{-1}
a(t)^2 \Bigr\}\nonumber
\\
&&\qquad\le\exp \biggl\{- \bigl(1+o(1) \bigr) \biggl(\frac{4}{4-\alpha}
\biggr)^{{(4-\alpha)}/ {2}} \biggl(\frac{2}{\alpha
}
\biggr)^{\alpha/2} \frac{a(t)^{4-\alpha}}{2c(\gamma)\theta^2 \kappa(d,\alpha)} \biggr\}
\\
&&\qquad\le\exp \bigl\{-(d+v)\log t \bigr\}
\nonumber
\end{eqnarray}
for some $v>0$, whenever $t$ is large and the constant
$u$ [appearing in (\ref{u-7})] satisfies
$ u>\theta^{{4}/{(4-\alpha)}}h(d,\alpha)$.

The asymptotic bound (\ref{u-11})
also holds in the setting of Theorem~\ref{intro-41}
by the same calculation of $\sigma_t^2$,
where (\ref{a-8}) in Lemma~\ref{a-6} is replaced
by (\ref{a-26}) in Lemma~\ref{a-24}.

By (\ref{u-10}), for large $t$ there is $v'>0$ such that
\[
\P \Bigl\{\max_{z\in2r\Z^d\cap Q_{2ta(t)+2r}}X_z(t)> 1 \Bigr\} \le\exp
\bigl\{-v'\log t \bigr\}.
\]
Consequently,
\[
\sum_k\P \Bigl\{\max_{z\in2r\Z^d\cap Q_{2t_ka(t_k)+2r}}X_z(t_k)>
1 \Bigr\} <\infty
\]
for $t_k=2^k$ ($k=1,2,\ldots $). By Borel--Cantelli lemma,
\[
\limsup_{k\to\infty}\max_{z\in2r\Z^d\cap
Q_{2t_ka(t_k)+2r}}X_z(t_k)
\le 1,\qquad  \mbox{a.s.}
\]
In view of (\ref{u-7}) and (\ref{u-9}),
\[
\limsup_{k\to\infty}(\log t_k)^{-{2}/{(4-\alpha)}}
\lambda_{\theta
V}(Q_{t_k}) \le \biggl(\frac{K}{ r}+1
\biggr)u,\qquad \mbox{a.s.}
\]
for any $u>\theta^{{4}/{(4-\alpha)}}h(d,\alpha)$.
Thus, (\ref{bound-21})
follows from the facst
that $\lambda_{\theta V}(Q_{t})$
is monotonic in $t$, $K>0$ is independent of $r$,
$r$ can be arbitrarily large
and $u$ can be arbitrarily close to $\theta^{{4}/{(4-\alpha)}}h(d,\alpha)$.

Based on the same argument, to establish (\ref{bound-23}) all we need
is to
show that
%
%
%e2.38 #&#
%
\begin{equation}
\label{u-12} \P \Bigl\{\sup_{g\in{\cal G}_1(Q_{r+1})} \bigl\langle V,
g_t^2\bigr\rangle\ge\theta^{-1}a(t)^2
\Bigr\} \le\exp \bigl\{- (1+v)\log t \bigr\}
\end{equation}
for some $v>0$, whenever $t$ is large and and
$ u>\frac{1}{2} (\frac{3}{2} )^{2/3}\theta^{4/3}$.

Indeed,
\[
\sigma_t^2=\sup_{g\in{\cal G}_1(Q_{r+1})}\int
_{-(r+1)}^{r+1} \bigl(g_t^2(x)
\bigr)^2\,dx \le a(t)\sup_{g\in{\cal G}_1(\R)}\int
_{-\infty}^{\infty}g^4(x)\,dx =\frac{3}{4}
\biggl(\frac{1}{2} \biggr)^{3/2}a(t),
\]
where the last step follows from (\ref{a-29}) in Lemma~\ref{a-27}.
By (\ref{u-10'}), therefore,
\begin{eqnarray*}&&\P \Bigl\{\sup_{g\in{\cal G}_1(Q_{r+1})}\bigl\langle
V, g_t^2\bigr\rangle\ge \theta ^{-1}a(t)^2
\Bigr\} \\
&&\qquad\le\exp \biggl\{- \bigl(1+o(1) \bigr) \biggl(\frac{2}{3}
\biggr)2^{3/2} \theta^{-2}a(t)^3 \biggr\}
\\
&&\qquad=\exp \biggl\{- \bigl(1+o(1) \bigr) \biggl(\frac{2}{3}
\biggr)2^{3/2}\theta ^{-2}u^{3/2} \log t \biggr\},
\end{eqnarray*}
which leads to (\ref{u-12}). \end{pf}

\begin{remark*}Clearly, (\ref{bound-21}) and
(\ref{bound-23}) still hold when $\lambda_{\theta V}(Q_t)$
is replaced by $\lambda_{\theta V}^+(Q_t)$.
Further, they can be improved into equalities where the limsup
can be strengthened into limit. The needed lower bounds will be given
in Lemma~\ref{l-3} below.
\end{remark*}

%s3 #&#
\section{Upper bounds}\label{u}

In this section we establish the upper bounds needed for
Theorems~\ref{intro-33}, \ref{intro-41} and~\ref{intro-46}.
Thanks to the homogeneity of the potential, the distribution of
the quenched moment in our theorems does not depends on $B_0$.
Therefore, we may take $B_0=0$ in the proof.
In other words, we prove that for any $\theta>0$,
%
%
%e3.1 #&#
%
\begin{eqnarray}
\label{u-1} &&\limsup_{t\to\infty}t^{-1}(\log
t)^{-{2}/{(4-\alpha)}}\log\E _0\exp \biggl\{ \theta\int
_0^tV(B_s)\,ds \biggr\}
\nonumber
\\[-8pt]
\\[-8pt]
\nonumber
&&\qquad\le
\theta^{{4}/ {(4-\alpha)}}h(d, \alpha),\qquad\mbox{a.s.}
\end{eqnarray}
in the context of Theorems~\ref{intro-33} or~\ref{intro-41},
where $h(d,\alpha)$ is defined in
(\ref{bound-22})
and
%
%
%e3.2 #&#
%
\begin{eqnarray}
\label{u-1'}&& \limsup_{t\to\infty}t^{-1}(\log
t)^{-2/3}\log\E_0\exp \biggl\{ \theta\int
_0^t V(B_s)\,ds \biggr\}
\nonumber
\\[-8pt]
\\[-8pt]
\nonumber
&&\qquad\le
\frac{1}{2} \biggl(\frac{3}{2} \biggr)^{2/3}
\theta^{4/3},\qquad \mbox{a.s.}
\end{eqnarray}
in the context of Theorem~\ref{intro-46}.

First, in all settings,
%
%
%e3.3 #&#
%
\begin{equation}
\label{exp-4} \E\otimes\E_0\exp \biggl\{\theta\int
_0^tV(B_s)\,ds \biggr\}<\infty,  \qquad t>0.
\end{equation}
Consequently,
\[
\E_0\exp \biggl\{\theta\int_0^tV(B_s)\,ds
\biggr\}<\infty, \qquad\mbox{a.s. } t>0.
\]
Here we recall our notation that ``$\E$,'' ``$\P$'' are used
for the expectation
and probability with respect to the Gaussian potential, and that
``$\E_0$,'' ``$\P_0$'' are used for the expectation
and probability with respect to the Brownian motion starting at 0.

Indeed, by the (conditional) Gaussian property stated in Lemma~\ref{def},
\[
\E\otimes\E_0\exp \biggl\{\theta\int_0^tV(B_s)\,ds
\biggr\} =\E_0\exp \biggl\{\frac{\theta^2}{2}\int
_0^t \int_0^t
\gamma (B_u-B_v)\,du\,dv \biggr\}.
\]
Therefore,
(\ref{exp-4}) follows from Theorem~4.3, \cite{CM-1} in the setting
of Theorem~\ref{intro-33}; from~(\ref{a-17}) below in the setting of
Theorem~\ref{intro-41}; and from Theorem~4.2.1, page~103, \cite{Chen}
in the setting of Theorem~\ref{intro-46}.

For any
open domain $D\in\R^d$, set the exit time
\[
\tau_D=\inf\{s\ge0; B_s\notin D\}.
\]
Recall the notation $Q_R=(-R, R)^d$.

In light of Lemma~\ref{bound-20}, our strategy for both upper and
lower bounds can be roughly outlined by
the following asymptotic relation:
%
%
%e3.4 #&#
%
\begin{equation}
\label{semi} \E_0\exp \biggl\{\theta\int_0^t
V(B_s)\,ds \biggr\} \approx\exp \bigl\{t\lambda_{\theta V}(Q_{R(t)})
\bigr\},
\end{equation}
where the principal
eigenvalue is introduced in (\ref{bound-16}), and square radius
$R(t)$ is nearly linear and carefully chosen according to the
context. To implement the upper bound, we consider the decomposition
\begin{eqnarray*} &&\E_0\exp \biggl\{\theta\int
_0^tV(B_s)\,ds \biggr\}
\\
&&\qquad=\E_0 \biggl[\exp \biggl\{\theta\int_0^t
V(B_s)\,ds \biggr\}; \tau_{Q_{R_1}}\ge t \biggr]
\\
&&\qquad\quad{}+\sum_{k=1}^\infty\E_0 \biggl[
\exp \biggl\{\theta\int_0^t V(B_s)\,ds
\biggr\}; \tau_{Q_{R_k}}< t\le\tau_{Q_{R_{k+1}}} \biggr]
\nonumber
\\
&&\qquad\le\E_0 \biggl[\exp \biggl\{\theta\int_0^t
V(B_s)\,ds \biggr\}; \tau_{Q_{R_1}}\ge t \biggr]
\\
&&\qquad\quad{}+\sum_{k=1}^\infty \bigl(\P_0
\{\tau_{Q_{R_k}}< t\} \bigr)^{1/2} \biggl\{\E_0 \biggl[
\exp \biggl\{2\theta\int_0^tV(B_s)\,ds
\biggr\}; \tau_{Q_{R_{k+1}}}\ge t \biggr] \biggr\}^{1/2},
\end{eqnarray*}
where
\begin{eqnarray*} R_k=\cases{ %
 \bigl(Mt (\log t)^{{1}/ {(4-\alpha)}} \bigr)^k,\vspace*{2pt}\cr
 \qquad $\mbox{in the context of Theorems~\ref{intro-33} or~\ref{intro-41},}$
\vspace*{2pt}\cr
\bigl(Mt (\log t)^{1/3} \bigr)^k,\vspace*{2pt}\cr
\qquad $\mbox{in the context of
Theorem~\ref{intro-46},}$} \qquad k=1,2,\ldots
\end{eqnarray*}
and the constant $M>0$ is fixed (for a while at least), but arbitrary.

The first term in the above decomposition is the dominating term
and is estimated in the following.
Let $p, q>1$ with $p^{-1}+q^{-1}=1$ with $p$ close to 1.
By
Lemma~4.3 [(4.5), with $\delta=1$ and $(\alpha, \beta)$
being replaced by $(p, q)$] in \cite{C}, we have for any $\varepsilon>0$,
\begin{eqnarray*} &&\E_0 \biggl[\exp \biggl\{\theta\int
_0^t V_\varepsilon(B_s)\,ds
\biggr\}; \tau_{Q_{R_1}}\ge t \biggr]
\\
&&\qquad\le \biggl(\E_0\exp \biggl\{\theta q\int_0^1V_\varepsilon(B_s)\,ds
\biggr\} \biggr)^{1/q}
\\
&&\qquad\quad{}\times \biggl\{\frac{1}{(2\pi)^{d/2}}\int
_{Q_{R_1}}\E_x \biggl[\exp \biggl\{ p\theta \int
_0^{t-1} V_\varepsilon(B_s)\,ds
\biggr\}; \tau_{Q_{R_1}}\ge t-1 \biggr]\,dx \biggr\}^{1/p},
\end{eqnarray*}
where the Gaussian field $V_\varepsilon(\cdot)$ is defined in
(\ref{app-3}).

The purpose of taking the above steps is to localize the Brownian range
and to re-shuffle the starting point
of the Brownian motion uniformly over $Q_{R_1}$. The Brownian motion
reaches anywhere of a super-linear (in $t$) distance from the origin with
a super-exponentially small probability which is negligible in comparison
to the essentially linear deviation scales shown in our main theorems.
The reason behind re-shuffling
is the explicit bounds (see, e.g., Lemmas~4.1 and 4.2 in \cite{C})
between the principal eigenvalues appearing in Lemma~\ref{bound-20}
and the exponential moment of the Brownian occupation time, in the case
when the Brownian motion has a uniformly distributed starting point.
Indeed, according to Lemma~4.1 in \cite{C},
\begin{eqnarray*}
&&\int_{Q_{R_1}}\E_x \biggl[\exp \biggl\{p\theta\int
_0^{t-1} V_\varepsilon(B_s)\,ds
\biggr\}; \tau_{Q_{R_1}}\ge t-1 \biggr]\,dx\\
&&\qquad \le\vert Q_{R_1}\vert
\exp \bigl\{(t-1)\lambda_{p\theta V_\varepsilon
}(Q_{R_1}) \bigr\}.
\end{eqnarray*}
Hence,
\begin{eqnarray*} &&\E_0 \biggl[\exp \biggl\{\theta\int
_0^t V_\varepsilon(B_s)\,ds
\biggr\}; \tau_{Q_{R_1}}\ge t \biggr]\,dx
\\
&&\qquad\le \biggl(\frac{2R_1^2}{\pi} \biggr)^{{d}/ {(2p)}} \biggl(
\E_0\exp \biggl\{ q\theta \int_0^1
V_\varepsilon(B_s)\,ds \biggr\} \biggr)^{1/q} \exp \bigl
\{(t-1)\lambda_{\theta pV_\varepsilon}(Q_{R_1}) \bigr\}.
\end{eqnarray*}

The reason for considering $V_\varepsilon$ instead of $V$ is that
Lemmas~4.3 and~4.1 in \cite{C}
were designed only for the pointwise defined functions.
To pass the above inequality
from $V_\varepsilon$ to $V$, we let $\varepsilon\to0^+$ on the
both sides. First notice that for any fixed $t$, by
comparing the variance between $V_\varepsilon$ and $V$, we have that
\begin{eqnarray*}
&&\E\otimes\E_0\exp \biggl\{\theta\int_0^t
V_\varepsilon(B_s)\,ds \biggr\}
\\
&&\qquad\le\E\otimes\E_0\exp
\biggl\{\theta\int_0^t V(B_s)\,ds
\biggr\}
\end{eqnarray*}
and by (\ref{exp-4}),
the right-hand side is finite for arbitrary $\theta>0$. Hence, a standard
argument by uniform integrability together with Lemma~\ref{def} leads to
%
%
%e3.5 #&#
%
\begin{equation}
\label{exp-5} \lim_{\varepsilon\to0^+}\E\otimes\E_0 \biggl\vert \exp
\biggl\{\theta\int_0^t V_\varepsilon(B_s)\,ds
\biggr\} -\exp \biggl\{\theta\int_0^t
V(B_s)\,ds \biggr\}\biggr \vert=0.
\end{equation}
Applying Fatou's lemma and
(\ref{lemma-2}) in Lemma~\ref{bound-14''}
to the inequality,
\begin{eqnarray}&&\E_0 \biggl[\exp \biggl\{\theta\int
_0^t V(B_s)\,ds \biggr\};
\tau_{Q_{R_1}}\ge t \biggr]\,dx
\nonumber\\
&&\qquad\le \biggl(\frac{2R_1^2}{\pi} \biggr)^{{d} /{(2p)}} \biggl(
\E_0\exp \biggl\{q\theta \int_0^1V(B_s)\,ds
\biggr\} \biggr)^{1/q}\exp \bigl\{(t-1) \lambda_{\theta pV}^+(Q_{R_1})
\bigr\},\nonumber\\
\eqntext{\mbox{a.s.}}
\end{eqnarray}
By a similar argument with $p=q=2$,
\begin{eqnarray} &&\E_0 \biggl[\exp \biggl\{2\theta\int
_0^tV(B_s)\,ds \biggr\};
\tau_{Q_{R_{k+1}}}\ge t \biggr]
\nonumber\\
&&\quad\le \biggl(\frac{2R_{k+1}^2}{\pi} \biggr)^{d/4} \biggl(\E_0
\exp \biggl\{4\theta\int_0^1
V(B_s)\,ds \biggr\} \biggr)^{1/2} \exp \bigl\{(t-1)
\lambda_{4\theta V}^+(Q_{R_{k+1}}) \bigr\},\nonumber\\
 \eqntext{\mbox{a.s.}}
\end{eqnarray}
for $k=1,2,\ldots.$

Summarizing our estimate,
\begin{eqnarray*}&&\E_0\exp \biggl\{\theta\int
_0^tV(B_s)\,ds \biggr\}
\\
&&\quad\le \biggl(\frac{2R_1^2}{\pi} \biggr)^{{d} /{(2p)}} \biggl(
\E_0\exp \biggl\{\theta q\int_0^1V(B_s)\,ds
\biggr\} \biggr)^{1/q} \exp \bigl\{(t-1)\lambda_{\theta pV}^+(Q_{R_1})
\bigr\}
\\
&&\quad\quad{}+ \biggl(\E_0\exp \biggl\{4\theta\int_0^1
V(B_s)\,ds \biggr\} \biggr)^{1/2}
\\
&&\phantom{\mbox{$+$}}\qquad{}\times \sum_{k=1}^\infty \biggl(
\frac{2R_{k+1}^2}{\pi} \biggr)^{d/4} \bigl(\P_0\{
\tau_{Q_{R_k}}< t\} \bigr)^{1/2}\exp \bigl\{(t-1)
\lambda_{4\theta V}^+(Q_{R_{k+1}}) \bigr\},\qquad\hspace*{-7pt} \mbox{a.s.}
\end{eqnarray*}

By the classic fact on the Gaussian tail,
\[
\bigl(\P_0\{\tau_{Q_{R_k}}< t\} \bigr)^{1/2} \le\exp
\bigl\{-c R_k^2/t \bigr\} =\exp \bigl\{-c
M^{2k}t^{2k-1}(\log t)^{{2k}/ {(4-\alpha)}} \bigr\}.
\]
Consequently, (\ref{u-1}) and (\ref{u-1'}) follow from
Lemma~\ref{bound-20}. Indeed, by (\ref{bound-21}) or~(\ref{bound-23}) (depending on the context),
the second term
(in the form of
infinite series) on the right-hand side of the established bound
is almost surely bounded when $M$ is sufficiently large,
and the first term contributes essentially up to the bound
given in (\ref{u-1}) or (\ref{u-1'})
as $p>1$ can be made arbitrarily close to~1.

%s4 #&#
\section{Lower bounds}\label{l}

In this section we establish the lower bounds needed for
Theorems~\ref{intro-33}, \ref{intro-41} and~\ref{intro-46}.
In other words, we prove that for any $\theta>0$,
%
%
%e4.1 #&#
%
\begin{eqnarray}
\label{l-1} &&\liminf_{t\to\infty}t^{-1}(\log
t)^{-{2}/{(4-\alpha)}}\log\E _0\exp \biggl\{ \theta\int
_0^tV(B_s)\,ds \biggr\}
\nonumber
\\[-8pt]
\\[-8pt]
\nonumber
&&\quad\ge
\theta^{{4}/ {(4-\alpha)}}h(d, \alpha), \qquad \mbox{a.s.}
\end{eqnarray}
in the context of Theorems~\ref{intro-33} or~\ref{intro-41},
where $h(d,\alpha)$ is defined in
(\ref{bound-22})
and
%
%
%e4.2 #&#
%
\begin{eqnarray}
\label{l-1'} &&\liminf_{t\to\infty}t^{-1}(\log
t)^{-2/3}\log\E_0\exp \biggl\{\int_0^t
V(B_s)\,ds \biggr\}
\nonumber
\\[-8pt]
\\[-8pt]
\nonumber
&&\quad\ge\frac{1}{2} \biggl(\frac{3}{2}
\biggr)^{2/3}\theta^{4/3},\qquad \mbox{a.s.}
\end{eqnarray}
in the context of Theorem~\ref{intro-46}.

Our treatment consists of two parts: Implementation of (\ref{semi})
for its lower bounds and establishment of the lower bounds for
the principal eigenvalues
which correspond to the upper bounds given in Lemma~\ref{bound-20}.

All notation used in Sections~\ref{bound} and~\ref{u} is adopted here.
Let $p, q>1$
satisfy $p^{-1}+q^{-1}=1$ with $p$ being close to 1, and let
$0<b<1$ be close to 1. For each $\varepsilon>0$, let
the pointwise defined potential
$V_\varepsilon(x)$ be given as
$(\ref{app-3})$.
Taking $\alpha=p$ and $q=\beta$, $\delta=t^b$ in
(4.6), Lemma~4.3, \cite{C} we have
\begin{eqnarray*}&&\E_0\exp \biggl\{\theta\int
_0^tV_\varepsilon(B_s)\,ds
\biggr\}
\\
&&\qquad\ge \biggl(\E_0\exp \biggl\{-\frac{q}{ p}\theta\int
_0^{t^b}V_\varepsilon (B_s)\,ds
\biggr\} \biggr)^{-p/q}\\
&&\qquad\quad{}\times \biggl( \int_{Q_{t^b}}p_{t^{b}}(x)
\E_x\exp \biggl\{\frac{\theta}{ p} \int_0^{t-t^b}V_\varepsilon(B_s)
\,dx \biggl\} \biggr)^p
\\
&&\qquad\ge \biggl(\E_0\exp \biggl\{-\frac{q}{ p}\theta\int
_0^{t^b}V_\varepsilon (B_s)\,ds
\biggr\} \biggr)^{-p/q} \\
&&\qquad\quad{}\times\biggl(\frac{e^{-ct^b}}{(2\pi t^b)^{d/2}} \int
_{Q_{t^b}}\E_x\exp \biggl\{\frac{\theta}{ p} \int
_0^{t-t^b}V_\varepsilon(B_s) \biggr
\}\,dx \biggr)^p,
\end{eqnarray*}
where $p_{t^b}(x)$ is the probability density of $B_{t^b}$.

Taking $\delta=t^b$ again
and replacing $t$, $\alpha$ and $\beta$ by $t-t^b$,
$p$ and $q$, respectively, in Lemma~4.2, \cite{C},
\begin{eqnarray*}&&\int_{Q_{t^b}}\E_x
\exp \biggl\{\frac{\theta}{ p} \int_0^{t-t^b}V_\varepsilon(B_s)
\biggr\}\,dx
\\
&&\qquad\ge(2\pi)^{pd/2}\bigl(t-t^b\bigr)^{{db} /{2}}
\bigl(t-t^b\bigr)^{{pd} /{(2q)}}\bigl(t-t^b
\bigr)^{-2db}
\\
&&\qquad\quad{}\times \exp \biggl\{-\frac{p}{ q}\bigl(t-t^b
\bigr)^b\lambda_{(p/q)\theta V_\varepsilon}(Q_{t^b}) \biggr\}\exp \bigl
\{pt\lambda_{\theta V_\varepsilon/p}(Q_{t^b}) \bigr\}.
 \end{eqnarray*}
Noticing that $\lambda_{\theta V_\varepsilon/p}(Q_{t^b}),
\lambda_{(p/q)\theta V_\varepsilon}(Q_{t^b})\ge0$, and replacing
$e^{-ct^b}$ by $e^{-Ct^b}$ for a larger $C$ to absorb all bounded-by-polynomial
quantities,
\begin{eqnarray*} \E_0\exp \biggl\{\theta\int
_0^tV_\varepsilon(B_s)\,ds
\biggr\}& \ge& e^{-Ct^b} \biggl(\E_0\exp \biggl\{-
\frac{q}{ p}\theta\int_0^{t^b}V_\varepsilon
(B_s)\,ds \biggr\} \biggr)^{-p/q}
\\
&&{}\times \exp \biggl\{-\frac{p^2}{ q}t^b \lambda_{(p/q)\theta V_\varepsilon}(Q_{t^b})
\biggr\}\exp \bigl\{t\lambda_{\theta V_\varepsilon/p}(Q_{t^b}) \bigr\}.
\end{eqnarray*}

Letting $\varepsilon\to0^+$ and taking the relation $V\stackrel{ d}{=}-V$
into account,
by (\ref{exp-5}) and~(\ref{lemma-2}) in
Lemma~\ref{bound-14''},
\begin{eqnarray*} &&\E_0\exp \biggl\{\theta\int
_0^tV(B_s)\,ds \biggr\}\\
&&\qquad \ge
e^{-Ct^b} \biggl(\E_0\exp \biggl\{-\frac{q}{ p}\theta
\int_0^{t^b}V (B_s)\,ds \biggr\}
\biggr)^{-p/q}
\\
&&\qquad\quad{}\times \exp \biggl\{-\frac{p^2}{ q}t^b \lambda_{(p/q)\theta V}^+(Q_{t^b})
\biggr\}\exp \bigl\{t\lambda_{\theta V/p}(Q_{t^b}) \bigr\}, \qquad\mbox{a.s.}
\end{eqnarray*}

Here we try to explain the strategy used in the above steps.
The Brownian motion is allowed to re-shuffle
its starting point uniformly over $Q_{t^b}$
within the affordable price $e^{-C t^b}$.
We take $b<1$ to make sure that
the energy
spent by the Brownian motion
during the ``relocation period'' $[0,t^b]$
is insignificant. Indeed, replacing $V$ by $-V$ and $t$ by $t^b$
in (\ref{u-1}) or in (\ref{u-1'}),
\[
\log\E_0\exp \biggl\{-\frac{q}{ p}\theta\int
_0^{t^b}V (B_s)\,ds \biggr\}=o(t),\qquad \mbox{a.s. }
(t\to\infty).
\]
In addition, by Lemma~\ref{bound-0},
\[
\frac{p^2}{ q}t^b
\lambda_{(p/q)\theta V}^+(Q_{t^b}) =o(t), \qquad\mbox{a.s.}
\]
under $b<1$.

On the other hand, we make $b$ close to 1 to
give the Brownian motion a decent chance to reach any location (within
the period $[0, t^b]$) up to
the distance $t^b\approx t$ where
the energy is rich to the degree requested by the lower bounds
in (\ref{l-1}) and~(\ref{l-1'}).

By the fact that
$p>1$ and $b<1$
can be made arbitrarily close to 1
[In particular, $\lambda_{\theta V/p}(Q_{t^b})\approx\lambda_{\theta
V}(Q_{t})$.],
the lower bounds
(\ref{l-1}) and (\ref{l-1'}) follow from the next lemma
which states another side of the story stated in
Lemma~\ref{bound-20}. %\end{pf}

%
%le4.1 #&#
\begin{lemma}\label{l-3} Under the assumptions of Theorems~\ref{intro-33}
or~\ref{intro-41}, for any $\theta>0$,
%
%
%e4.3 #&#
%
\begin{equation}
\label{l-3'} \liminf_{t\to\infty}(\log
t)^{-{2}/{(4-\alpha)}}\lambda_{\theta V}(Q_t) \ge
\theta^{{4}/ {(4-\alpha)}}h(d,\alpha),\qquad \mbox{a.s.},
\end{equation}
where $h(d,\alpha)$ is given in (\ref{bound-22}).

Under the assumption of Theorem~\ref{intro-46}, for any $\theta>0$
%
%
%e4.4 #&#
%
\begin{equation}
\label{l-3''} \liminf_{t\to\infty}(\log
t)^{-2/3}\lambda_{\theta V} \bigl((-t, t) \bigr) \ge
\frac{1}{2} \biggl(\frac{3}{2} \biggr)^{2/3}
\theta^{4/3},\qquad \mbox{a.s.}
\end{equation}
\end{lemma}

\begin{pf} Recall that
$a(t)$ and $g_t(x)$ are defined in (\ref{u-7}) and (\ref{u-8}),
respectively. Let the constant $r>0$ be fixed but
arbitrary, and set
${\cal N}_t=2r\Z^d\cap Q_{t-r}$.
By~(\ref{u-7'}) and by the monotonicity of $\lambda_{\theta V}(D)$
in the set $D\subset\R^d$,
\begin{eqnarray*}\lambda_{\theta V}(Q_t)\ge
a(t)^2\max_{z\in{\cal N}_t} \sup_{g\in{\cal F}_d( a(t)z+Q_r)}
\biggl\{\theta a(t)^{-2}\bigl\langle V, g_t^2
\bigr\rangle -\frac{1}{2} \int_{a(t)z+Q_r}\bigl\vert\nabla g(x)
\bigr\vert^2\,dx \biggr\}.
 \end{eqnarray*}

For any $g\in{\cal G}_d(Q_r)$ and $z\in{\cal N}_t$, notice that
$g^z(\cdot)\equiv g(\cdot-a(t)z)\in{\cal F}_d(a(t)z+Q_r)$,
and by translation invariance,
\[
\int_{a(t)z+Q_r}\bigl\vert\nabla g^z(x)
\bigr\vert^2\,dx =\int_{Q_r}\bigl\vert\nabla g(x)
\bigr\vert^2\,dx,\qquad z\in{\cal N}_t.
\]
Consequently,
%
%
%e4.5 #&#
%
\begin{equation}
\label{l-5} \lambda_{\theta V}(Q_t)\ge a(t)^2
\biggl\{ \theta a(t)^{-2}\max_{z\in{\cal N}_t}\bigl\langle V,
\bigl(g^z\bigr)_t^2\bigr\rangle-
\frac{1}{2} \int_{Q_r}\bigl\vert\nabla g(x)
\bigr\vert^2\,dx \biggr\}
\end{equation}
for any $g\in{\cal F}_d(Q_r)$.
In the following argument $g\in{\cal F}_d(Q_r)$ is fixed but arbitrary.
Set $t_k=2^k$ ($k=1,2,\ldots $). Our next step
is to show that
%
%
%e4.6 #&#
%
\begin{equation}
\label{l-7'} \liminf_{k\to\infty}a(t_k)^{-2}
\max_{z\in{\cal N}_{t_k}}\bigl\langle V, \bigl(g^z
\bigr)_{t_k}^2\bigr\rangle \ge\sigma(g), \qquad\mbox{a.s.}
\end{equation}
whenever
%
%
%e4.7 #&#
%
\begin{equation}
\label{l-7} \cases{ %
u< \bigl(2dc(\gamma)
\bigr)^{{2} /{(4-\alpha)}},
&\quad $\mbox{in the context of Theorem~\ref{intro-33}},$
\vspace*{2pt}\cr
u< (2d C_H )^{{2}/ {(4-\alpha)}},&\quad $\mbox{in the context of
Theorem~\ref{intro-41}},$
\vspace*{2pt}\cr
u<2^{2/3},&\quad $\mbox{in the context of Theorem~\ref{intro-46}},$
}
\end{equation}
where
\[
\sigma(g)=\cases{ %
\displaystyle\biggl( \int_{\R^d\times\R^d}\frac{g^2(x)g^2(y)}{\vert x- y\vert
^p}\,dx\,dy \biggr)^{1/2},\vspace*{2pt}\cr\qquad
  $\mbox{in Theorem~\ref{intro-33}},$
\vspace*{2pt}\cr
\displaystyle\Biggl(\int_{\R^d\times\R^d}g^2(x)g^2(y)
\Biggl(\prod_{j=1}^d\vert
x_j-y_j\vert^{2-2H_j} \Biggr)^{-1}\,dx\,dy
\Biggr)^{1/2},\vspace*{2pt}\cr\qquad $\mbox{in Theorem~\ref{intro-41}},$
\vspace*{2pt}\cr
\displaystyle\biggl(\int_{-\infty}^\infty g^4(x)\,dx
\biggr)^{1/2},\vspace*{2pt}\cr\qquad $\mbox{in Theorem~\ref{intro-46}.}$}
\]

The proof of (\ref{l-7'}) in the setting of Theorem~\ref{intro-46} is
easy due to the fact that the sequence
\[
\bigl\langle V, \bigl(g^z\bigr)_t^2\bigr
\rangle,\qquad z\in{\cal N}_t
\]
is an i.i.d. family with the common distribution $N (0,
a(t)\sigma^2(g) )$.
Consequently,
\begin{eqnarray*}&&\P \biggl\{\max_{z\in{\cal N}_t}\bigl
\langle V, \bigl(g^z\bigr)_t^2\bigr\rangle
\le a(t)^2 \biggl(\int_{-\infty}^\infty
g^4(x)\,dx \biggr)^{1/2} \biggr\}
\\
&&\qquad= \biggl(1-\P \biggl\{\bigl\langle V, \bigl(g^0
\bigr)_t^2\bigr\rangle>a(t)^2 \biggl(\int
_{-\infty}^\infty g^4(x)\,dx
\biggr)^{1/2} \biggr\} \biggr)^{\#({\cal N}_t)}.
 \end{eqnarray*}
By the classic tail estimate for normal distribution,
\begin{eqnarray*}&&\P \biggl\{\bigl\langle V, \bigl(g^0
\bigr)_t^2\bigr\rangle >a(t)^2 \biggl(\int
_{-\infty}^\infty g^4(x)\,dx
\biggr)^{1/2} \biggr\}
\\
&&\qquad=\exp \biggl\{- \bigl(1+o(1) \bigr)\frac{a(t)^3}{2} \biggr\} =\exp \biggl\{-
\bigl(1+o(1) \bigr)\frac{u^{3/2}\log t}{2} \biggr\}.
\end{eqnarray*}
By the fact that $\#({\cal N}_t)\sim(2r)^{-1}t$ as $t\to\infty$, we have
%
%
%e4.8 #&#
%
\begin{equation}
\label{l-5'} \P \biggl\{\max_{z\in{\cal N}_t}\bigl\langle
V, \bigl(g^z\bigr)_t^2\bigr\rangle \le
a(t)^2 \biggl(\int_{-\infty}^\infty
g^4(x)\,dx \biggr)^{1/2} \biggr\}\le\exp \bigl
\{-t^\beta \bigr\}
\end{equation}
for some $\beta>0$, whenever $u <2^{2/3}$. Consequently,
\[
\sum_k\P \biggl\{\max_{z\in{\cal N}_{t_k}}
\bigl\langle V, \bigl(g^z\bigr)_{t_k}^2\bigr
\rangle \le a(t_k)^2 \biggl(\int_{-\infty}^\infty
g^4(x)\,dx \biggr)^{1/2} \biggr\}<\infty.
\]
Hence, (\ref{l-7'}) follows from Borel--Cantelli lemma.

In the settings of Theorems \ref{intro-33} and~\ref{intro-41},
the proof of (\ref{l-7'}) is harder due to lack of independence.
Our approach relies on the control of the covariance. Write
$\xi_z(t)=\langle V, (g^z)_t^2\rangle$.
For each $z,z'\in{\cal N}_t$,
\begin{eqnarray*}
&&\Cov(\xi_z, \xi_{z'})\\
&&\qquad=\int
_{\R^d\times\R^d}\gamma(x-y) \bigl(g^z\bigr)_t^2(x)
\bigl(g^{z'}\bigr)_t^2(y)\,dx\,dy
\\
&&\qquad=\int_{\R^d\times\R^d}\gamma \bigl(x-y+\bigl(z-z'\bigr)
\bigr)g_t^2(x)g_t^2(y)\,dx\,dy
\\
&&\qquad=\int_{\R^d\times\R^d}\gamma \bigl(a(t)^{-1}(x-y)+
\bigl(z-z'\bigr) \bigr)g^2(x)g^2(y)\,dx\,dy, \qquad z,
z'\in{\cal N}_t.
\end{eqnarray*}

Taking $z=z'$ in the setting of Theorem~\ref{intro-33},
%
%
%e4.9 #&#
%
\begin{eqnarray}
\label{l-8} \Var \bigl(\xi_0(t) \bigr)&=&\int_{\R^d\times\R^d}
\gamma \bigl(a(t)^{-1}(x-y) \bigr) g^2(x)g^2(y)\,dx\,dy
\nonumber
\\[-8pt]
\\[-8pt]
\nonumber
&\sim& c(\gamma)\sigma^2(g)a(t)^\alpha\qquad (t\to\infty),
\end{eqnarray}
where the last step follows from (\ref{intro-15}).\vadjust{\goodbreak}

Using (\ref{intro-39}) instead of (\ref{intro-15}),
we can see that in the setting of Theorem~\ref{intro-41},
%
%
%e4.10 #&#
%
\begin{equation}
\label{l-9} \Var \bigl(\xi_0(t) \bigr)= C_H
\sigma^2(g)a(t)^\alpha\qquad (t>0).
\end{equation}
We now claim that in both settings,
%
%
%e4.11 #&#
%
\begin{equation}
\label{l-10} R_t\equiv\mathop{\max_{ z, z'\in{\cal
N}_t}}_{
z\not=z'}
\bigl\vert\Cov \bigl(\xi_z(t), \xi_{z'}(t) \bigr)\bigr \vert =o
\bigl(a(t)^\alpha \bigr)\qquad (t\to\infty).
\end{equation}

By the assumption that $\gamma(x)$ is bounded on $\{\vert x\vert\ge
1\}$,
$\Cov (\xi_z(t),\xi_{z'}(t) )$ is bounded uniformly over the pairs
$(z, z')$ with $z\not=z'$ and over $t$ in the setting of Theorem~\ref{intro-33}. In particular, (\ref{l-10}) holds.

The proof of (\ref{l-10})
is a little trickier when it comes to Theorem~\ref{intro-41}. That is
the reason why we cannot have a constant bound
for $\Cov (\xi_z(t),\xi_{z'}(t) )$ with $z\not=z'$.
More precisely, $\Cov (\xi_z(t),\xi_{z'}(t) )\to\infty$
as $t\to\infty$ when $z_j=z_j'$ for some $1\le j\le d$. Here
we use the notation $z=(z_1,\ldots, z_d)$.
Write
\[
J\bigl(z, z'\bigr)=\bigl\{1\le j\le d; z_j=z_j'
\bigr\},\qquad  z, z'\in{\cal N}_t.
\]
By (\ref{intro-39}),
\begin{eqnarray*}
&&\int_{\R^d\times\R^d}\gamma \bigl(a(t)^{-1}(x-y)+
\bigl(z-z'\bigr) \bigr)g^2(x)g^2(y)\,dx\,dy
\\
&&\qquad\sim C_H \biggl(\prod_{j\notin J(z,z')}\bigl\vert
z_j-z_j'\bigr\vert ^{2-2H_j}
\biggr)^{-1} a(t)^{\alpha(z,z')}
\\
&&\qquad\quad{}\times\int_{\R^d\times\R^d}g^2(x)g^2(y)
\biggl(\prod_{j\in J(z,z')}\vert x_j-y_j
\vert^{2-2H_j} \biggr)^{-1}\,dx\,dy\qquad (t\to\infty),
\end{eqnarray*}
where
\[
\alpha\bigl(z, z'\bigr)=\sum_{j\in J(z,z')}(2-2H_j).
\]
By the fact that $\vert z_j-z_j'\vert\ge2r$ for $j\notin J(z, z')$,
the above asymptotic equivalence can be developed into the uniform bound
\[
\mathop{\max_{ z, z'\in{\cal N}_t}}_{
z\not=z'} \bigl\vert\Cov \bigl(
\xi_z(t), \xi_{z'}(t) \bigr) \bigr\vert \le C
a(t)^{\alpha'},
\]
where
\[
\alpha'\equiv\mathop{\max_{ z, z'\in{\cal
N}_t}}_{z\not=z'}
\alpha\bigl(z, z'\bigr)<\alpha.
\]
So (\ref{l-10}) holds.

Given a small but fixed $v>0$, taking $A=\sigma(g)a(t)^2$ and
$B=v\sigma(g) a(t)^2$ in
Lemma~\ref{l-14} below,
\begin{eqnarray*}
&&\P \Bigl\{\max_{z\in{\cal N}_t}\xi_z(t)\le\sigma(g)
a(t)^2 \Bigr\}
\\
&&\qquad\le \biggl(\P \biggl\{\xi_0(t)\le(1+v)\sigma(g) a(t)^2
\sqrt\frac{2R_t+\Var (\xi_0(t) )}{\Var (\xi
_0(t) )} \biggr\} \biggr)^{\#{\cal N}_t}\\
&&\qquad\quad{} +\P \bigl\{U\ge v\sigma(g)
a(t)^2/\sqrt{2 R_t} \bigr\},
\end{eqnarray*}
where $U$ is a standard normal random variable.

For the second term on the right-hand side,
\[
\P \bigl\{U\ge v\sigma(g) a(t)^2 \bigr\} =\exp \biggl\{-
\bigl(1+o(1) \bigr)\frac{v^2a(t)^4\sigma^2(g)}{4R_t} \biggr\} \le\exp \{-2\log t \}
\]
for large $t$, where the last step follows from (\ref{l-10}).

As for the first term, by (\ref{l-8}) and (\ref{l-9}) the algorithm used
in (\ref{l-5'}) shows that it is bounded by $e^{-t^\beta}$
for some $\beta>0$ when $t$ is large, $v$ is small and $u$
satisfies (\ref{l-7}).

Summarizing our computation, we obtain a bound that leads to
\[
\sum_k\P \Bigl\{\max_{z\in{\cal N}_{t_k}}
\xi_z(t_k)\le\sigma(g) a(t_k)^2
\Bigr\} <\infty.
\]
So (\ref{l-7'}) follows from
Borel--Cantelli
lemma.

In view of (\ref{l-5}), (\ref{l-7'}) implies that for every
$g\in{\cal F}_d(Q_r)$,
%
%
%e4.12 #&#
%
\begin{eqnarray}
\label{l-11} &&\liminf_{k\to\infty}(\log t_k)^{-{2}/{(4-\alpha)}}
\lambda _{\theta
V}(Q_{t_k})
\nonumber
\\
&&\qquad\ge \bigl(2d c(\gamma) \bigr)^{{2}/ {(4-\alpha)}} \biggl\{ \theta
\biggl(\int_{\R^d\times\R^d}\frac{g^2(x)g^2(y)}{\vert
x-y\vert
^\alpha}\,dx\,dy \biggr)^{1/2}\\
&&\hspace*{163pt}{}-
\frac{1}{2}\int_{\R^d}\bigl\vert\nabla g(x)
\bigr\vert^2\,dx \biggr\},\qquad  \mbox{a.s.}\nonumber
\end{eqnarray}
in the setting of Theorem~\ref{intro-33}, that
%
%
%e4.13 #&#
%
\begin{eqnarray}
\label{l-12} &&\hspace*{-4pt}\liminf_{k\to\infty}(\log t_k)^{-{2}/{(4-\alpha)}}
\lambda _{\theta
V}(Q_{t_k})
\nonumber\\
&&\hspace*{-7pt}\qquad\ge(2dC_H)^{{2}/ {(4-\alpha)}} \Biggl\{ \theta \Biggl(
\int_{\R^d\times\R^d}g^2(x)g^2(y) \Biggl(\prod
_{j=1}^d\vert x_j-y_j
\vert^{2-2H_j} \Biggr)^{-1}\,dx\,dy \Biggr)^{1/2}
\nonumber
\\[-8pt]
\\[-8pt]
\nonumber
&&\hspace*{256pt}\hspace*{-4pt}{} -\frac{1}{2}\int_{\R^d}\bigl\vert\nabla g(x)
\bigr\vert^2\,dx \Biggr\},\\
 \eqntext{\mbox{a.s.}}
\end{eqnarray}
in the setting of Theorem~\ref{intro-41}, and that
%
%
%e4.14 #&#
%
\begin{eqnarray}
\label{l-13} &&\liminf_{k\to\infty}(\log t_k)^{-2/3}
\lambda_{\theta V}(Q_{t_k})
\nonumber
\\[-8pt]
\\[-8pt]
\nonumber
&&\qquad\ge2^{2/3} \biggl\{\theta
\biggl(\int_{-\infty}^\infty g^4(x)\,dx
\biggr)^{1/2} -\frac{1}{2}\int_{-\infty}^\infty
\bigl\vert g'(x)\bigr\vert^2\,dx \biggr\},\qquad  \mbox{a.s.}
\end{eqnarray}
in the setting of Theorem~\ref{intro-46}.

By the monotonicity of $\lambda_{\theta V}(Q_t)$ in $t$, the liminf
along the
sub-sequence $t_k$ in~(\ref{l-11}), (\ref{l-12}) and (\ref{l-13})
can be extended into the liminf along the continuous time $t$.

Recall that $W^{1,2}(\R^d)$ is the Sobolev space
defined in (\ref{intro-31'}). Consistently with~(\ref{bound-0}), we define
\[
{\cal F}_d\bigl(\R^d\bigr)= \bigl\{g\in W^{1,2}
\bigl(\R^d\bigr); \|g\|_2=1 \bigr\}.
\]
We now prove that the functions $g$ on the right-hand sides
of (\ref{l-11}), (\ref{l-12}) and (\ref{l-13})
can be extended from ${\cal F}_d(Q_r)$ to ${\cal F}_d(\R^d)$,
and complete the proof of
Lemma~\ref{l-3}.

We start with (\ref{l-11}).
The right-hand side
can be extended to all $g\in{\cal F}_d(\R^d)$ for
the following two reasons: First, the infinitely smooth, rapidly decreasing
and locally supported functions are dense in the Sobolev space
$W^{1,2}(\R^d)$ under the Sobolev norm
\[
\|g\|_{W^{1,2}(\R^d)}\equiv\sqrt{\|g\|_2^2+
\tfrac{1}{2}\|\nabla g\|_2^2}
\]
and $r>0$ in (\ref{l-11}) is arbitrary.
Second, by (\ref{a-2}) the functional
\[
{\cal F}(g)= \biggl( \int_{\R^d\times\R^d}\frac{g^2(x)g^2(y)}{\vert x-y\vert^\alpha}\,dx\,dy
\biggr)^{1/2}-\frac{1}{2}\int_{\R^d}\bigl\vert\nabla
g(x)\bigr\vert^2\,dx
\]
is continuous under the Sobolev norm $\|\cdot\|_{W^{1,2}(\R^d)}$.

Taking supremum over $g\in{\cal F}_d(\R^d)$ on the right-hand side of
(\ref{l-11}) we obtain the lower bound
\begin{eqnarray*}
&&\liminf_{t\to\infty}(\log t)^{-{2}/{(4-\alpha)}}\lambda_{\theta V}(Q_{t})
\\
&&\qquad\ge \bigl(2d c(\gamma) \bigr)^{{2}/ {(4-\alpha)}}M_{d,\alpha}(
\theta)
\\
&&\qquad=\frac{4-\alpha}{4} \biggl(\frac{\alpha}{2} \biggr)^{{\alpha} /{(4-
\alpha)}} \bigl(2dc(\gamma)\theta^2\kappa(d,\alpha)
\bigr)^{{2}/ {(4-\alpha)}},\qquad \mbox{a.s.}
\end{eqnarray*}
in the setting of Theorem~\ref{intro-33},
where $M_{d,\alpha}(\theta)$ is defined in (\ref{a-4}),
and the last step follows from the variation identity (\ref{a-7}).

Using (\ref{a-25}) (with $\alpha_j=2-2H_j$) instead of (\ref{a-7}),
by the same argument, from (\ref{l-12}) we derive that
\begin{eqnarray*}
&&\liminf_{t\to\infty}(\log t)^{-{2}/{(4-\alpha)}}\lambda_{\theta V}(Q_{t})
\\
&&\qquad\ge (2dC_H )^{{2} /{(4-\alpha)}}\widetilde{M}_{d,\alpha
}(
\theta)
\\
&&\qquad=\frac{4-\alpha}{4} \biggl(\frac{\alpha}{2} \biggr)^{{\alpha} /{(4-
\alpha)}} \bigl(2dc(\gamma)\theta^2\tilde{\kappa}(d, H)
\bigr)^{{2} /{( 4-\alpha)}}, \qquad \mbox{a.s.}
\end{eqnarray*}
in the setting of Theorem~\ref{intro-41}.

In the same way, by (\ref{a-28}) and (\ref{l-13}) we have
\begin{eqnarray*}
&&\liminf_{t\to\infty}(\log t)^{2/3}\lambda_{\theta V}(Q_{t})
\\
&&\qquad\ge2^{2/3}\sup_{g\in{\cal F}_1(\R)} \biggl\{\theta \biggl( \int
_{-\infty}^\infty g^4(x)\,dx
\biggr)^{1/2} -\frac{1}{2}\int_{-\infty}^\infty
\bigl\vert g'(x)\bigr\vert^2\,dx \biggr\}
\\
&&\qquad=\frac{1}{2} \biggl(\frac{3}{2} \biggr)^{2/3}
\theta^{4/3}, \qquad \mbox{a.s.}
\end{eqnarray*}
in the setting of Theorem~\ref{intro-46}.
\end{pf}

We end this section with the following lemma.

%
%le4.2 #&#
\begin{lemma}\label{l-14}
Let $(\xi_1,\ldots, \xi_n)$ be a mean-zero Gaussian vector with
identically distributed components.
Write
\[
R=\max_{i\not=j}\bigl\vert\Cov (\xi_i,
\xi_j)\bigr\vert
\]
and assume that $\Var(\xi_1)\ge2R$.
Then for any $A, B>0$,
\[
\P \Bigl\{\max_{k\le n}\xi_k\le A \Bigr\} \le
\biggl(\P \biggl\{\xi_1 \le\sqrt\frac{ 2R+\Var(\xi_1)}{\Var(\xi_1)}(A +B) \biggr\}
\biggr)^n +\P\{U\ge B/\sqrt{2R}\},
\]
where $U$ is a standard normal random variable.
\end{lemma}

\begin{pf} Let $\eta_1,\ldots, \eta_n$ be an i.i.d. sequence
independent of
$U$. Assume that
$\eta_1\stackrel{ d}{=}\xi_1$ and write
\[
\zeta_k=\sqrt\frac{ \Var(\xi_1)}{2R+\Var(\xi_1)}(\eta_k+\sqrt{2R}U).
\]
With the assumption $\Var(\xi_1)\ge2R$,
it is straightforward to exam that
\[
\Var(\xi_k)=\Var(\zeta_k) \quad\mbox{and}\quad \Cov(
\xi_i, \xi_j)\le\Cov(\zeta_i,
\zeta_j),\qquad  i,j, k=1,\ldots, n.
\]
By Slepian's lemma (\cite{Slepian}, see also Lemma~5.5.1, \cite{MR}),
\[
\P \Bigl\{\max_{k\le n}\xi_k\le A \Bigr\}\le \P
\Bigl\{\max_{k\le n}\zeta_k\le A \Bigr\}.
\]
Notice that
\[
\max_{k\le n}\zeta_k=\sqrt\frac{ 2R\Var(\xi_1)}{2R+\Var(\xi_1)}U+
\sqrt\frac{\Var(\xi_1)}{2R+\Var(\xi_1)}\max_{k\le n}\eta_k.
\]
By the triangle inequality,
\[
\P \Bigl\{\max_{k\le n}\xi_k\le A \Bigr\}\le\P
\biggl\{ \max_{k\le n}\eta_k\le\sqrt
\frac{ 2R+\Var(\xi_1)}{\Var(\xi
_1)}(A+B) \biggr\} +\P \{U\le-B/\sqrt{2R}\}.
\]
The conclusion follows from the symmetry of $U$ and
the independence of~$\{\eta_k\}$.
\end{pf}

\begin{appendix}\label{app}
%s5 #&#
\section*{Appendix}

%s5.1 #&#
\subsection{Brownian integral as a limit}%\label{app}

In this subsection, $\langle V, \varphi\rangle$ ($\varphi\in{\cal
S}(\R^d)$)
is a mean-zero generalized Gaussian field with homogeneity defined
in (\ref{intro-9}). Let $\mu(dx)$ be the spectral measure of
$\langle V, \varphi\rangle$ and let the pointwise defined Gaussian field
$V_\varepsilon(x)$ ($x\in\R^d$) be given in (\ref{app-3}). The main
goal here is to
prove
%
%le5.1 #&#
\setcounter{lemmas}{0}
\begin{lemmas}\label{def} Assume that
%
%
%e5.1 #&#
%
\setcounter{equation}{0}
\begin{equation}
\label{app-0} \int_{\R^d}\frac{1}{1+ \vert\lambda\vert^2}\mu(d\lambda)<
\infty.
\end{equation}
Under the product law $\P\otimes\P_x$, the ${\cal L}^2$-limit
%
%
%e5.2 #&#
%
\begin{equation}
\label{def-1} \int_0^tV(B_s)\,ds
\stackrel{\mathit{def}} { =}\lim_{\varepsilon\to0^+}\int_0^tV_\varepsilon(B_s)\,ds
\end{equation}
exists for every $t\ge0$. In addition, there is a modification of
the limiting process in~(\ref{def-1}) that is
$(\frac{1}{2}-u)$-H\"older continuous for any $u>0$.
Further,
conditioned on the Brownian motion, the process
%
%
%e5.3 #&#
%
\begin{equation}
\label{app-1} \int_0^tV(B_s)\,ds,\qquad t
\ge0
\end{equation}
is mean-zero Gaussian with the (conditional) variance
%
%
%e5.4 #&#
%
\begin{equation}
\label{def-2}\quad  \E \biggl\{\int_0^tV(B_s)\,ds
\biggr\}^2= \frac{1}{(2\pi)^d} \int_{\R^d} \biggl\vert
\int_0^t e^{i\lambda\cdot B_s}\,ds \biggr\vert ^2
\mu (d\lambda),\qquad  t\ge0.
\end{equation}
\end{lemmas}

\begin{pf} First notice that conditioned on the Brownian motion, the process
\[
I_\varepsilon(t)=\int_0^tV_\varepsilon(B_s)\,ds,\qquad
t\ge0
\]
is Gaussian with the conditional variance
%
%
%e5.5 #&#
%
\begin{equation}
\label{app-5} \E I_\varepsilon^2(t) =\frac{1}{(2\pi)^d}\int
_{\R^d} \biggl\vert\int_0^t
e^{i\lambda\cdot B(u)}\,du \biggr\vert^2\bigl\vert{\cal F}(l) (\varepsilon \lambda )
\bigr\vert^2 \mu(d\lambda).
\end{equation}

We claim that there is a constant $C>0$ such that
%
%
%e5.6 #&#
%
\begin{equation}
\label{app-6} \int_{\R^d}\E_x \biggl\vert\int
_0^t e^{i\lambda\cdot B(u)}\,du \biggr\vert^2\mu(d
\lambda) \le C\bigl(t\vee t^2\bigr),\qquad t\ge0.
\end{equation}
Indeed,
\begin{eqnarray*}
\E_x \biggl\vert\int_0^t
e^{i\lambda\cdot B(u)}\,du \biggr\vert^2 &=&\int_0^t
\int_0^t\E_x e^{i\lambda\cdot(B_u-B_v)}\,du\,dv
\\
&=&\int_0^t \int_0^t
\exp \biggl\{-\frac{\vert\lambda\vert^2}{2}\vert u-v\vert \biggr\}\,du\,dv.
\end{eqnarray*}
The right-hand side is equal to
\[
\frac{4}{\vert\lambda\vert^2} \biggl[t-\frac{2}{\vert\lambda\vert^2} \bigl(1-e^{-t\vert\lambda\vert^2/2} \bigr)
\biggr],
\]
which yields a bound $4t/\vert\lambda\vert^2$ for $\vert\lambda
\vert\ge1$.
As for $\vert\lambda\vert\le1$, we use the trivial bound
\[
\int_0^t\!\! \int_0^t
\exp \biggl\{-\frac{\vert\lambda\vert^2}{2}\vert u-v\vert \biggr\}\,du\,dv \le t^2.
\]

Thus
\[
\int_{\R^d}\E_x \biggl\vert\int_0^t
e^{i\lambda\cdot B(u)}\,du \biggr\vert^2\mu(d\lambda) \le4t\int
_{\{\vert\lambda\vert\ge1\}}\frac{1}{\vert\lambda\vert
^2}\mu (d\lambda) +t^2\int
_{\{\vert\lambda\vert\le1\}}\mu(d\lambda).
\]
Hence, (\ref{app-6}) follows from (\ref{app-0}).

To prove the ${\cal L}^2$-convergence described in (\ref{def-1}),
all we need is to establish the
existence of the limit
$ \lim_{\varepsilon,\varepsilon'\to0^+}
\E_x\otimes\E (I_{\varepsilon'} (t)I_\varepsilon(t) )$.

Indeed, similar to (\ref{app-5}),
\[
\E_x\otimes\E \bigl(I_{\varepsilon'} (t)I_\varepsilon(t)
\bigr) =\frac{1}{(2\pi)^d}\int_{\R^d}\E_x\biggl \vert
\int_0^t e^{i\lambda\cdot B(u)}\,du \biggr\vert^2
{\cal F}(l) (\varepsilon\lambda)\ol{{\cal F}(l) \bigl(\varepsilon'\lambda
\bigr)}\mu(d\lambda).
\]
By (\ref{app-6}), the fact that
\[
\bigl\vert{\cal F}(l) (\varepsilon\lambda)\bigr\vert\le1 \quad\mbox{and}\quad \lim_{\varepsilon\to0^+}{
\cal F}(l) (\varepsilon \lambda)=1,
\]
and by the dominant convergence theorem we obtain
\[
\lim_{\varepsilon,\varepsilon'\to0^+} \E_x\otimes\E \bigl(I_{\varepsilon'}
(t)I_\varepsilon(t) \bigr) =\frac{1}{(2\pi)^d}\int_{\R^d}
\E_x\biggl \vert\int_0^t
e^{i\lambda\cdot B(u)}\,du \biggr\vert^2 \mu(d\lambda).
\]

Write $I_0(t)= \lim_{\varepsilon\to0^+}I_\varepsilon(t)$
as the ${\cal L}^2(\P_x\otimes\P)$-limit. Recall the classical
fact that the ${\cal L}^2$-limit
of Gaussian process remains Gaussian. Conditioned on the Brownian
motion,
$\{I_0(t);  t\ge0\}$ is Gaussian with zero mean and
the conditional variance
%
%
%e5.7 #&#
%
\begin{equation}
\label{app-7} \E I_0^2(t)=\frac{1}{(2\pi)^d}\int
_{\R^d} \biggl\vert\int_0^t
e^{i\lambda\cdot B(u)}\,du \biggr\vert^2 \mu(d\lambda).
\end{equation}

Strictly speaking, $\{I_0(t);  t\ge0\}$ exists
as a family of equivalent classes. In the following we try to find
a continuous modification of this family. For any $s,t\ge0$ with $s<t$,
notice that $I_0(t)- I_0(s)$ is conditionally normal with the variance
\begin{eqnarray*}
\E \bigl[I_0(t)- I_0(s) \bigr]^2&=&
\frac{1}{(2\pi)^d}\int_{\R^d} \biggl\vert\int_s^te^{i\lambda\cdot
B(u)}\,du
\biggr\vert^2 \mu(d\lambda)\\
&\stackrel{ d} {=}&\frac{1}{(2\pi)^d}\int
_{\R^d}\biggl \vert\int_0^{t-s}e^{i\lambda\cdot B(u)}\,du
\biggr\vert^2\mu (d\lambda).
\end{eqnarray*}
Thus, for any integer $m\ge1$,
%
%
%e5.8 #&#
%
\begin{eqnarray}
\label{app-8}&& \E_x\otimes\E \bigl[I_0(t)-
I_0(s) \bigr]^{2m}
\nonumber
\\[-8pt]
\\[-8pt]
\nonumber
&&\qquad= (2m-1)!! \E_x \biggl(\int
_{\R^d}\biggl \vert\int_0^{t-s}e^{i\lambda\cdot B(u)}
\,du \biggr\vert^2\mu(d\lambda) \biggr)^m.
\end{eqnarray}

To estimate the right-hand side, we consider the nonnegative,
continuous process
\[
Z_t= \biggl\{\int_{\R^d} \biggl\vert\int
_0^{t}e^{i\lambda\cdot B(u)} \,du \biggr\vert^2
\mu(d\lambda) \biggr\}^{1/2},\qquad t\ge0.
\]
By the triangle inequality,
%
%
%e5.9 #&#
%
\begin{equation}
\label{app-8'} Z_{s+t}\le Z_t+Z_s',\qquad
s,t\ge0,
\end{equation}
where
\[
Z_s'= \biggl\{\int_{\R^d} \biggl\vert
\int_t^{t+s}e^{i\lambda\cdot B(u)} \,du\biggr \vert^2
\mu(d\lambda) \biggr\}^{1/2}
\]
is independent of $\{B_u;  0\le u\le t\}$ and equal in law to
$Z_s$. By (1.3.7), page 21 in~\cite{Chen}, for any $t, a, b> 0$,
\[
\P_x\{Z_t\ge a+b\}\le\P_x
\{Z_t\ge a\}\P_x\{Z_t\ge b\}.
\]
Consequently,
\[
\P_x\{Z_t\ge Mn\sqrt{t} \}\le \bigl(\P_x
\{Z_t\ge M\sqrt{t} \} \bigr)^n,\qquad  n=1,2,\ldots.
\]
By (\ref{app-6}), one can take $M>0$ sufficiently large so
\[
\sup_{0<t\le1}\P_x\{Z_t\ge M\sqrt{t}\}
\le e^{-2}.
\]
Hence,
%
%
%e5.10 #&#
%
\begin{equation}
\label{app-9} \sup_{0<t\le1}\E_x\exp \bigl
\{M^{-1}Z_t/\sqrt{t} \bigr\}<\infty.
\end{equation}
Replacing $t$ by $t-s$ and applying it to (\ref{app-8}), we obtain
\[
\E\otimes\E\bigl\vert I_0(t)-I_0(s)\bigr\vert^{2m}\le
C_m\vert t-s\vert ^m\qquad \mbox{for all $s, t\ge0$ with $
\vert t-s\vert\le1$}.
\]
By the classic result on chaining (see, e.g., Lemma~9, \cite{CLR}),
there is a modification
of $\{I_0(t);  t\ge0\}$ that is $(\frac{1}{2}-u)$-H\"older continuous
for any $u>0$. \end{pf}

%
%le5.2 #&#
\begin{lemmas}\label{intro-25'} Under the assumptions in Theorems~\ref
{intro-33}, \ref{intro-41} or~\ref{intro-46}, (\ref{intro-50}) holds
for some $\delta>0$.
In particular, the Brownian integral in (\ref{app-1}) is well-defined as
stated in Lemma~\ref{def}.
\end{lemmas}

\begin{pf} We first consider the setting of
Theorem~\ref{intro-33}. By the fact that $\mu$ is tempered,
all we need to show is
\[
\int_{\{\vert\lambda\vert\ge1\}} \frac{1}{\vert\lambda\vert^{2(1-\delta)}}\mu(d\lambda)<\infty.
\]

Let $\varphi$ be the density of the standard normal distribution on
$\R^d$. By Fourier transform
\begin{eqnarray*}
2^{kd}\int_{\R^d}\gamma(x)\varphi
\bigl(2^kx\bigr)\,dx&=&\frac{1}{(2\pi)^d} \int_{\R^d}
\exp \biggl\{-\frac{\vert2^{-k}\lambda\vert^2}{2} \biggr\}
\mu(d\lambda)\\
& \ge& c\mu \bigl
\{2^{k-1}\le\vert\lambda\vert\le2^k \bigr\}.
\end{eqnarray*}
On the other hand, by (\ref{intro-15})
\begin{eqnarray*}
2^{kd}\int_{\R^d}\gamma(x)\varphi
\bigl(2^kx\bigr)\,dx&=&\int_{\R^d}\gamma
\bigl(2^{-k}x\bigr)\varphi(x)\,dx\\
& \sim& c(\gamma)2^{\alpha k}\int
_{\R^d}\frac{\varphi(x)}{\vert x\vert
^\alpha}\,dx\qquad (k\to\infty).
\end{eqnarray*}
Hence, there is a constant $C>0$ such that
\[
\mu \bigl\{2^{k-1}\le\vert\lambda\vert\le2^k \bigr\}\le C
2^{\alpha k},\qquad k=1,2,\ldots.
\]
Thus
\[
\int_{\{\vert\lambda\vert\ge1\}}\frac{1}{\vert\lambda\vert
^{2(1-\delta
)}}\mu(d\lambda) \le C\sum
_{k=1}^\infty2^{-2(1-\delta)(k-1)} \mu \bigl
\{2^{k-1}\le\vert\lambda\vert\le2^k \bigr\}<\infty
\]
for any $\delta< \frac{2-\alpha}{2}$.

In the setting of Theorem~\ref{intro-46}, where
$\mu(d\lambda)=d\lambda$ is the 1-dimensional Lebesgue measure,
the validity of (\ref{intro-50})
can be directly verified with any $\delta<1$.

As for the setting of Theorem~\ref{intro-41}, by (\ref{intro-39})
and spherical substitution,
\begin{eqnarray*} \int_{\R^d}\frac{1}{(1+ \vert\lambda\vert^2)^{1-\delta}}\mu
(d\lambda)& =&\widehat{C}_H\int_{\R^d} \Biggl(\prod
_{j=1}^d\vert\lambda_j\vert
^{2H_j-1} \Biggr)^{-1} \frac{1}{(1+ \vert\lambda\vert^2)^{1-\delta}}\,d\lambda
\\
&=&C\int_0^\infty r^{-(d-\alpha)}
\frac{r^{d-1}}{(1+ r^2)^{1-\delta}}\,dr\\
& =&C\int_0^\infty
\frac{r^{\alpha-1}}{(1+ r^2)^{1-\delta}}\,dr<\infty
\end{eqnarray*}
as $\delta<\frac{2-\alpha}{2}$. \end{pf}

%s5.2 #&#
\subsection{Counting the covering balls}

Let $D, D'\subset\R^d$ be two domains in $\R^d$
and ${\cal Q}(D)$ be a class of functions on $D$. Assume that $D'$
is bounded. For
each $\varepsilon>0$, let $\rho_\varepsilon(f, g)$ be a pseudometric on
${\cal Q}(D)$
such that
\[
\rho_\varepsilon(f,g)\le \biggl(\int_{D'}\bigl\vert
A_\varepsilon (f) (x)-A_\varepsilon(g) (x) \bigr\vert \,dx \biggr)^{1/2}
\Bigl(\sup_{x\in D'} \bigl\vert B_\varepsilon(f) (x)-
B_\varepsilon(g) (x)\bigr\vert \Bigr)^{1/2},
\]
where $A_\varepsilon$ and $B_\varepsilon$ are two (possibly nonlinear) maps
from ${\cal Q}(D)$ to the space $\mbox{Lip}(D')$ of Lipschitz functions
on $D'$. Assume further that there are constants $C>0$, $p>1$, $m\ge1$
such that
%
%
%e5.11 #&#
%
\begin{eqnarray}
\label{ball-1}\bigl \vert B_\varepsilon(g) (x)\bigr\vert&\le& C
\quad\mbox{and}
\nonumber
\\[-8pt]
\\[-8pt]
\nonumber
\bigl \vert
B_\varepsilon(g) (x)- B_\varepsilon(g) (y)\bigr\vert&\le& C\varepsilon^{-m}
\vert x-y\vert,\qquad x, y\in D',
\\
\label{ball-2} \int_{D'}\bigl\vert A_\varepsilon(g) (x)
\bigr\vert^p\,dx&\le& C\quad \mbox{and}
\nonumber
\\[-8pt]
\\[-8pt]
\nonumber
\bigl\vert A_\varepsilon(g) (x)-
A_\varepsilon(g) (y)\bigr\vert&\le& C\varepsilon^{-m}\vert x-y\vert,\qquad x, y\in
D'
\end{eqnarray}
uniformly for all $g\in{\cal Q}(D)$ and sufficiently small $\varepsilon>0$.

%
%le5.3 #&#
\begin{lemmas}\label{ball-3}
Under the above assumptions,
\[
\log N \bigl({\cal Q}(D),\rho_\varepsilon, \varepsilon \bigr)=O \biggl(
\varepsilon^{-{2p}/{(2p-1)}}\log\frac{1}{\varepsilon} \biggr) \qquad\bigl(\varepsilon\to0^+\bigr).
\]
\end{lemmas}

\begin{pf} Notice that
\begin{eqnarray*}
\hat{\rho}_\varepsilon(f,g)&=&\int_{D'}\bigl\vert
A_\varepsilon (f) (x)-A_\varepsilon(g) (x) \bigr\vert \,dx\quad
\mbox{and}\\
\rho_\varepsilon^*(f, g)&=&\sup_{x\in D'}\bigl\vert
B_\varepsilon(f) (x)- B_\varepsilon (g) (x)\bigr\vert
\end{eqnarray*}
define two pseudometrics on ${\cal Q}(D)$.
We now claim that for any $u, v>0$ with $\sqrt{uv}=\varepsilon$,
%
%
%e5.13 #&#
%
\begin{equation}
\label{ball-5} N \bigl({\cal Q}(D), \rho_\varepsilon, \varepsilon \bigr) \le N
\bigl({\cal Q}(D), \rho^*_\varepsilon, u \bigr) N \bigl({\cal Q}(D), \hat{
\rho}_\varepsilon, v \bigr).
\end{equation}

Indeed, we first cover ${\cal Q}(D)$ by
$N ({\cal Q}(D), \rho_\varepsilon^*, u )$
$\rho_\varepsilon^*$-balls with the diameter
smaller than $u$. For each such ball, it can be covered
by at most $N ({\cal Q}(D), \hat{\rho}_\varepsilon, v )$ of
$\hat{\rho}_\varepsilon$-balls with the diameter smaller than
$v$.
In this way, the set ${\cal Q}(D)$ is covered by at most
$N ({\cal Q}(D), \rho^*, u )
N ({\cal Q}(D), \hat{\rho}_\varepsilon, v )$ of its nonempty
subsets. For $f, g$ coming from same subset,
$\hat{\rho}_\varepsilon(f,g)< u$
and\break  $\rho_\varepsilon^*(f, g)< v$. Hence
\[
\rho_\varepsilon(f,g)\le\sqrt{\hat{\rho}_\varepsilon(f,g)\rho^*(f, g)} <
\sqrt{uv}=\varepsilon.
\]
Hence, (\ref{ball-5}) holds.

With (\ref{ball-5}), it is sufficient to establish
%
%
%e5.14 #&#
%
\begin{eqnarray}
\label{ball-6} &&N \bigl({\cal Q}(D), \rho_\varepsilon^*,
\varepsilon^{{2p}/
{(2p-1)} }\bigr)
\nonumber
\\[-8pt]
\\[-8pt]
\nonumber
&&\qquad=\exp \biggl\{O \biggl(\varepsilon^{-{2p}/{(2p-1)}}\log
\frac{1}{\varepsilon
} \biggr) \biggr\}\qquad \bigl(\varepsilon\to0^+\bigr)
\\
\label{ball-7}&& N \bigl({\cal Q}(D), \hat{\rho}_\varepsilon,
\varepsilon^{{(2(p-1))} /{ (2p-1)}} \bigr)
\nonumber
\\[-8pt]
\\[-8pt]
\nonumber
&&\qquad =\exp \biggl\{O \biggl(
\varepsilon^{-{2p}/{(2p-1)}}\log\frac{1}{\varepsilon
} \biggr) \biggr\} \qquad \bigl(\varepsilon
\to0^+\bigr).
\end{eqnarray}
Indeed, applying (\ref{ball-5}) with
%
%
%e5.16 #&#
%
\begin{equation}
\label{ball-8} u(\varepsilon) =\varepsilon^{{2p} /{(2p-1)} }\quad\mbox{and}\quad v(
\varepsilon)=\varepsilon^{{(2(p-1))} /{(2p-1)}},
\end{equation}
and using (\ref{ball-6}) and (\ref{ball-7}) we have
\begin{eqnarray*}N \bigl({\cal Q}(D), \rho_\varepsilon, \varepsilon
\bigr) &\le& N \bigl({\cal Q}(D), \rho_\varepsilon^*, u(\varepsilon) \bigr) N
\bigl({\cal Q}(D), \hat{\rho}_\varepsilon, v(\varepsilon) \bigr)
\\
&=&\exp \biggl\{O \biggl(\varepsilon^{-{2p}/{(2p-1)}}\log\frac{1}{\varepsilon
} \biggr)
\biggr\} \qquad\bigl(\varepsilon\to0^+\bigr).
\end{eqnarray*}

We first prove (\ref{ball-6}). Let $u(\varepsilon)$ be defined in (\ref
{ball-8}).
Define the map
\[
B_\varepsilon^*\dvtx {\cal Q}(D)\longrightarrow \bigl((\sqrt{d}C)^{-1}
\varepsilon^{m}u(\varepsilon) \Z^d\cap D'
\bigr)^{u(\varepsilon) \Z\cap[-C, C]}
\]
as $B_\varepsilon^* f(x)=u(\varepsilon)  [u(\varepsilon)^{-1}B_\varepsilon
(f)(x_0) ]$
whenever
\[
x\in \bigl(x_0-(2\sqrt{d}C)^{-1}\varepsilon^mu\bigl(
\varepsilon\bigr), x_0+(2\sqrt{d}C)^{-1}\varepsilon^mu(
\varepsilon) \bigr]^d
\]
for some
$x_0\in(\sqrt{d}C)^{-1}\varepsilon^mu(\varepsilon)\Z^d
\cap D'$, where $[\cdot]$ is the integer-part function.\vadjust{\goodbreak}

By (\ref{ball-1})
\[
\sup_{x\in D'}\bigl\vert B_\varepsilon g(x) -B_\varepsilon^*
g(x)\bigr\vert< \frac{u(\varepsilon)}{2},\qquad  g\in{\cal Q}(D).
\]
Consequently, for any $f, g\in{\cal Q}(D)$ with
$B_\varepsilon^* f=B_\varepsilon^* g$,
\[
\rho_\varepsilon^*(f,g)= \sup_{x\in D'}\bigl\vert
B_\varepsilon f(x)-B_\varepsilon g(x)\bigr\vert< u(\varepsilon).
\]
Hence,
\begin{eqnarray*} N \bigl({\cal Q}(D), \rho_\varepsilon^*, u(
\varepsilon) \bigr) &\le&\# \bigl\{ \bigl((\sqrt{d}C)^{-1}
\varepsilon^mu(\varepsilon) \Z^d\cap D'
\bigr)^{u(\varepsilon) \Z\cap[-C, C]} \bigr\}
\\
&=&\exp \biggl\{O \biggl(\varepsilon^{-{2p}/{(2p-1)}}\log\frac{1}{\varepsilon
} \biggr)
\biggr\}\qquad \bigl(\varepsilon\to0^+\bigr).
\end{eqnarray*}

It remains to establish (\ref{ball-7}).
Let $v(\varepsilon)$ be given in (\ref{ball-8}), and write
\[
M_\varepsilon= \bigl(8C v(\varepsilon)^{-1} \bigr)^{(p-1)^{-1}}.
\]
Define the map
\[
A^*_\varepsilon\dvtx  {\cal Q}(D)\longrightarrow \bigl(\bigl(4\bigl\vert
D'\bigr\vert \sqrt{d}C\bigr)^{-1}\varepsilon^mv(
\varepsilon) \Z^d\cap D' \bigr)^{(8\bigl\vert D'\bigr\vert)^{-1}v(\varepsilon) \Z\cap
[-M_\varepsilon,
M_\varepsilon]}
\]
as $A^*_{\varepsilon}g(x)
= \{(8\vert D'\vert)^{-1}v(\varepsilon) [8\vert D'\vert
v(\varepsilon
)^{-1} A_\varepsilon g(x_0)
]\wedge M_\varepsilon \}\vee(-M_\varepsilon)$,
whenever
\[
x\in \bigl(x_0-\bigl(8\bigl\vert D'\bigr\vert \sqrt{d}C
\bigr)^{-1}\varepsilon^mv(\varepsilon), x_0+\bigl(8
\bigl\vert D'\bigr\vert \sqrt{d}C\bigr)^{-1}\varepsilon^mv(
\varepsilon) \bigr]^d
\]
for some $x_0\in(4\vert D'\vert
\sqrt{d}C)^{-1}\varepsilon^mv(\varepsilon)
\Z^d\cap D'$.

By (\ref{ball-2}),
\begin{eqnarray*}&&\sup_{g\in{\cal Q}(D)}\int
_{D'} \bigl\vert A_\varepsilon(g) (x)-A_\varepsilon^*(g)
(x)\bigr \vert \,dx
\\
&&\qquad\le\frac{1}{4}v(\varepsilon)+2\sup_{g\in{\cal Q}(D)}
\int_{\{\vert A_\varepsilon(g)\vert>M_\varepsilon\}} \bigl\vert A_\varepsilon(g) (x)\bigr\vert \,dx
\\
&&\qquad\le\frac{1}{4}v(\varepsilon)+2M_\varepsilon^{-(p-1)}C\le
\frac{1}{
2}v(\varepsilon).
\end{eqnarray*}

Consequently, for $f, g\in{\cal Q}(D)$ with
$A_\varepsilon^*f=A_\varepsilon^*g$,
$\hat{\rho}_\varepsilon(f,g)<v(\varepsilon)$ for small $\varepsilon$.
Hence,
\begin{eqnarray*}N \bigl({\cal Q}(D), \hat{\rho}_\varepsilon, v(
\varepsilon) \bigr) &\le&\# \bigl\{ \bigl(\bigl(4\vert D'\vert \sqrt{d}C
\bigr)^{-1}\varepsilon^m v(\varepsilon) \Z^d\cap
D' \bigr)^{(8\vert D'\vert)^{-1}
v(\varepsilon) \Z\cap[-M_\varepsilon, M_\varepsilon]} \bigr\}
\\
&=&\exp \biggl\{O \biggl(\varepsilon^{-{2p}/{(2p-1)}}\log\frac{1}{\varepsilon
} \biggr)
\biggr\}\qquad \bigl(\varepsilon\to0^+\bigr).
\end{eqnarray*}
\upqed\end{pf}

%s5.3 #&#
\subsection{Variations}

In this section we establish some Sobolev-type inequalities
and validate the variations used in the paper. Recall
that $W^{1,2}(\R^d)$ is the Sobolev space defined in (\ref
{intro-31'}) and
\[
{\cal F}_d\bigl(\R^d\bigr)= \bigl\{g\in W^{1,2}
\bigl(\R^d\bigr); \|g\|_2=1 \bigr\}.
\]
Similar to (\ref{bound-1}), define
\[
{\cal G}_d\bigl(\R^d\bigr)= \bigl\{g\in
W^{1,2}\bigl(\R^d\bigr); \|g\|_2^2+
\tfrac{1}{2}\|g\|_2^2=1 \bigr\}.
\]

Recall
(Lemma~7.2, \cite{C}) that for any $0\le\alpha<2\wedge d$
there is $C_\alpha>0$ such that
%
%
%e5.17 #&#
%
\begin{equation}
\label{a-1} \int_{\R^d}\frac{f^2(x)}{\vert x\vert^\alpha}\,dx\le
C_\alpha\|f\| _2^{2-\alpha} \|\nabla f
\|_2^{\alpha},\qquad f\in W^{1,2}\bigl(\R^d\bigr).
\end{equation}
A simple trick by translation invariance, show that (\ref{a-1}) remains
true with
the same constant $C_\alpha$ if the left-hand side is replaced
by
\[
\sup_{y\in\R^d}\int_{\R^d}\frac{f^2(x)}{\vert x-y\vert^\alpha}\,dx.
\]
Immediately,
%
%
%e5.18 #&#
%
\begin{eqnarray}
\label{a-2} \int_{\R^d\times\R^d}\frac{f^2(x)f^2(y)}
{\vert x-y\vert^\alpha}\,dx\,dy& =&\int
_{\R^d}f^2(y) \biggl[\int_{\R^d}
\frac{f^2(x)}{\vert x-y\vert
^\alpha
}\,dx \biggr]\,dy
\nonumber
\\[-8pt]
\\[-8pt]
\nonumber
&\le& C_\alpha\|f\|_2^{4-\alpha}
\|\nabla f\|_2^{\alpha}
\end{eqnarray}
for every $f\in W^{1,2}(\R^d)$.

As a consequence, the constant
%
%
%e5.19 #&#
%
\begin{eqnarray}
\label{a-3} \kappa(d,\alpha)&=&\inf \biggl\{C>0; \int_{\R^d\times\R^d}
\frac{f^2(x)f^2(y)}{\vert x-y\vert^\alpha}\,dx\,dy \
\nonumber
\\[-8pt]
\\[-8pt]
\nonumber
&&\hspace*{20pt}\le C\|f\|_2^{4-\alpha} \|\nabla f
\|_2^{\alpha}\ \forall f\in W^{1,2}\bigl(
\R^d\bigr) \biggr\}
\end{eqnarray}
is finite.

Other variations relevant to Theorem~\ref{intro-33} are
%
%
%e5.20 #&#
%
\begin{eqnarray}
\label{a-4} M_{d,\alpha}(\theta)&=&\sup_{g\in{\cal F}_d(\R^d)} \biggl\{
\theta \biggl( \int_{\R^d\times\R^d}\frac{g^2(x)g^2(y)}{\vert x-y\vert^\alpha
}\,dx\,dy
\biggr)^{1/2}
\nonumber
\\[-8pt]
\\[-8pt]
\nonumber
&&\hspace*{96pt}{}-\frac{1}{2}\int_{\R^d}\bigl\vert\nabla
g(x)\bigr\vert^2\,dx \biggr\},\qquad \theta>0,
\\
\label{a-5}\qquad \sigma(d,\alpha)&=&\sup_{g\in{\cal G}_d(\R^d)} \biggl\{ \int
_{\R^d\times\R^d}\frac{g^2(x)g^2(y)}{\vert x-y\vert^\alpha
}\,dx\,dy \biggr\}^{1/2}.
\end{eqnarray}
By (\ref{a-1}), one can easily
show that $M_{d,\alpha}(\theta)$ and
$\sigma(d,\alpha)$ are finite under the assumption $0\le\alpha
<2\wedge d$.

%
%le5.4 #&#
\begin{lemmas}\label{a-6} Under $\alpha<2\wedge d$,
%
%
%e5.22 #&#
%
\begin{eqnarray}
\label{a-7} M_{d,\alpha}(\theta)&=&\frac{4-\alpha}{4} \biggl(\frac{\alpha}{2}
\biggr)^{{\alpha}/ {(4-\alpha)}} \kappa(d,\alpha)^{{2}/ {(4-
\alpha)}}\theta^{{4} /{(4-\alpha)}},
\\
\label{a-8} \sigma(d,\alpha)&=& \biggl(\frac{4-\alpha}{4} \biggr)^{{(4-
\alpha)} /{4}} \biggl(\frac{\alpha}{2} \biggr)^{\alpha/4} \kappa(d,
\alpha)^{1/2}.
\end{eqnarray}
\end{lemmas}

\begin{pf} Let $f\in{\cal F}(\R^d)$ be fixed but arbitrary, and let
$C_f>0$ satisfy
\[
\int_{\R^d\times\R^d}\frac{f^2(x)f^2(y)}{\vert x-y\vert^\alpha}\,dx\,dy=C_f \|\nabla
f\|_2^\alpha.
\]
Given $\beta>0$ let $g(x)=\beta^{d/2} f(\beta x)$. Then
$\|\nabla g\|_2=\beta\|\nabla f\|_2$ and therefore
\[
\int_{\R^d\times\R^d}\frac{g^2(x)g^2(y)}{\vert x-y\vert^\alpha}\,dx\,dy =\beta^\alpha\int
_{\R^d\times\R^d}\frac{f^2(x)f^2(y)}{\vert
x-y\vert
^\alpha}\,dx\,dy =C_f
\beta^\alpha\|\nabla f\|_2^\alpha.
\]
By the fact that
$g\in{\cal F}_d(\R^d)$,
\[
M_{d,\alpha}(\theta)\ge\theta C_f^{1/2}
\beta^{\alpha/2}\|\nabla f\|_2^{\alpha/2} -\tfrac{1}{2}\|
\nabla g\|_2^2 =\theta C_f^{1/2}
\beta^{\alpha/2}\|\nabla f\|_2^{\alpha/2} -\tfrac{1}{2}
\beta^2\|\nabla f\|_2^2.
\]
Notice $\beta\|\nabla f\|_2$ runs over all positive numbers. So we have
\[
M_{d,\alpha}(\theta)\ge\sup_{x>0} \biggl\{\theta
C_f^{1/2} x^{\alpha/2} -\frac{1}{2}x^2
\biggr\}=\frac{4-\alpha}{4} \biggl(\frac{\alpha}{2} \biggr)^{{
\alpha} /{(4-\alpha)}}C_f^{{2}/ { (4-\alpha)}}
\theta^{{4}/ {(4-\alpha)}}.
\]
Take supremum over $f$ on the right-hand side. Noticing
that ${\cal S}(\R^d)$ is dense in $W^{1,2}(\R^d)$, by space homogeneity
we have established the relation ``$\ge$'' for (\ref{a-7}).

On the other hand, for any $g\in{\cal F}_d(\R^d)$,
\begin{eqnarray*}&&\theta \biggl( \int_{\R^d\times\R^d}
\frac{g^2(x)g^2(y)}{\vert x-y\vert^\alpha
}\,dx\,dy \biggr)^{1/2} -\frac{1}{2}\int
_{\R^d}\bigl\vert\nabla g(x)\bigr\vert^2\,dx
\\
&&\qquad\le\theta\kappa(d,\alpha)^{1/2}\|\nabla g\|_2^{\alpha/2}-
\frac{1}{2} \|\nabla g\|_2^2 \le\sup
_{x>0} \biggl\{\theta\kappa(d,\alpha)^{1/2}
x^{\alpha/2} -\frac{1}{2}x^2 \biggr\}
\\
&&\qquad=\frac{4-\alpha}{4} \biggl(\frac{\alpha}{2} \biggr)^{{\alpha} /{(4-
\alpha)}} \kappa(d, \alpha)^{{2}/ {(4-\alpha)}}\theta^{{4}
/{(4-\alpha)}}.
\end{eqnarray*}
Taking supremum over $g\in{\cal F}_d(\R^d)$ on the left-hand side,
we reach the relation ``$\le$'' for (\ref{a-7}).

For any $g\in{\cal F}_d(\R^d)$ by space homogeneity,
\begin{eqnarray*}&&\frac{1}{\sigma(d,\alpha)} \biggl(\int_{\R^d\times\R^d}
\frac{g^2(x)g^2(y)
}{\vert x-y\vert^\alpha}\,dx\,dy \biggr)^{1/2} -\frac{1}{2}\int
_{\R^d}\bigl\vert\nabla g(x)\bigr\vert^2\,dx
\\
&&\qquad\le\frac{1}{\sigma(d,\alpha)}\sigma(d,\alpha) \bigl(1+\|\nabla g\|_2^2
\bigr)-\frac{1}{2}\|\nabla g\|_2^2=1.
\end{eqnarray*}
Taking supremum over $g$,
\[
M_{d,\alpha} \biggl(\frac{1}{\sigma(d,\alpha)} \biggr)\le1.
\]
Combining this with (\ref{a-7}) we have proved the ``$\ge$'' half
for (\ref{a-8}).

On the other hand, for any $f\in W^{1,2}(\R^d)$,
\begin{eqnarray*}&& \biggl(\int_{\R^d\times\R^d}
\frac{f^2(x)f^2(y)
}{\vert x-y\vert^\alpha}\,dx\,dy \biggr)^{1/2}\\
&&\qquad \le\kappa(d,\alpha)^{1/2}
\|f\|_2^{{(4-\alpha)}/ {2}}\|\nabla f\| _2^{\alpha/2}
\\
&&\qquad=\kappa(d,\alpha)^{1/2} \biggl(\frac{2\alpha}{4-\alpha}
\biggr)^{\alpha/4} \bigl(\|f\|_2^2
\bigr)^{{(4-\alpha)}/ {4}} \biggl(\frac{4-\alpha}{2\alpha}\|\nabla f
\|_2^2 \biggr)^{\alpha/4}
\\
&&\qquad\le\kappa(d,\alpha)^{1/2} \biggl(\frac{2\alpha}{4-\alpha}
\biggr)^{\alpha/4} \frac{4-\alpha}{4} \biggl(\|f\|_2^2+
\frac{1}{2}\|\nabla f\|_2^2 \biggr),
\end{eqnarray*}
where the last step follows from the H\"older inequality
$ab\le p^{-1}a^p+q^{-1}b^q$ with $p=4(4-\alpha)^{-1}$ and $q=4/\alpha$.
This leads to the ``$\le$'' half for (\ref{a-8}). \end{pf}

We need an inequality
comparable to the one in (\ref{a-1})
for formulating and proving
Theorem~\ref{intro-41}, but could not find it in literature.
We establish it in the following.

Let the real numbers $\alpha_1,\ldots, \alpha_d$ satisfy
$0\le\alpha_j<1$ and \mbox{$\alpha\equiv\alpha_1+\cdots+\alpha_d <2$}.

%
%le5.5 #&#
\begin{lemmas}\label{a-12} For any $\theta>0$,
%
%
%e5.24 #&#
%
\begin{equation}
\label{a-13} \qquad\sup_{g\in{\cal F}_d(\R^d)} \Biggl\{\theta\int
_{\R^d} \Biggl(\prod_{j=1}^d
\vert x_j\vert^{-\alpha_j} \Biggr) g^2(x)\,dx-
\frac{1}{2}\int_{\R^d}\bigl\vert \nabla g(x)
\bigr\vert^2\,dx \Biggr\}<\infty.
\end{equation}
\end{lemmas}

\begin{pf} Define the function
\[
K(x)=\prod_{j=1}^d\vert
x_j\vert^{-\alpha_j},\qquad  x=(x_1,\ldots,
x_d) \in\R^d.
\]
The fact that $K(x)$
blows up at every coordinate plane make the problem harder
comparing to setting of
the Newtonian kernel $\vert x\vert^{-\alpha}$ which blows
up only at $0$. The fact
that $\alpha_1,\ldots, \alpha_d$ are allowed to be different posts
an extra challenge.
The proof provided here is probabilistic.

Let the linear Brownian motions
$B_1(s),\ldots, B_d(s)$ be the independent components of the $d$-dimensional
Brownian motion $B_s$ and define the process
%
%
%e5.25 #&#
%
\begin{equation}
\label{a-14} \eta_t=\int_0^t
K(B_s)\,ds=\int_0^t \Biggl(\prod
_{j=1}^d\bigl\vert B_j(s)\bigr\vert
^{-\alpha_j} \Biggr)\,ds,\qquad t >0.
\end{equation}
This process is well defined under our assumption on $\alpha_1,\ldots
,\alpha_d$.
Indeed, it is not hard to see that for each $t>0$, $\E_0\eta_t<\infty$.
Further, we now prove that there is a $b>0$ such that
%
%
%e5.26 #&#
%
\begin{equation}
\label{a-15} \E_0\exp \bigl\{b\eta_1^{2/\alpha}
\bigr\}<\infty.
\end{equation}
We point out that (\ref{a-15}) is a strengthened version of the
exponential integrability for $\eta_1$ obtained by Hu, Nualart and Song
(Lemma A.5, \cite{HNS}) and the approach for~(\ref{a-15})
presented here is modified
from theirs.

Given the integer $m\ge1$,
\begin{eqnarray*}\E_0\eta_t^m&=&
\int_{[0,t]^m}\,ds_1\cdots ds_m\prod
_{j=1}^d\E_0\prod
_{k=1}^m \bigl\vert B_1(s_k)
\bigr\vert^{-\alpha_j}\\
&=& m!\int_{[0,t]_<^m}\,ds_1\cdots
ds_m\prod_{j=1}^d
\E_0\prod_{k=1}^m \bigl\vert
B_1(s_k)\bigr\vert^{-\alpha_j},
\end{eqnarray*}
where the multi-dimensional time set $[0,t]_<^m$ is defined as
\[
[0,t]_<^m= \bigl\{(s_1,\ldots s_m)
\in[0,t]^m; s_1<s_2\cdots <s_m
\bigr\}.
\]

Let $(s_1,\ldots, s_m)\in[0,t]_<^m$ be fixed for a while
and ${\cal A}_s=\sigma\{B_1(u);  0\le u\le s\}$ be the filtration
generated by the linear Brownian motion $B_1(t)$. Write
\begin{eqnarray*} \E_0 \bigl\{\bigl\vert B_1(s_k)
\bigr\vert^{-\alpha_j}|{\cal A}_{s_{k-1}} \bigr\} &=&\int_0^\infty
\P_0 \bigl\{\bigl\vert B_1(s_k)
\bigr\vert^{-\alpha_j}\ge a|{\cal A}_{s_{k-1}} \bigr\}\,da
\\
&=&\int_0^\infty\P_0 \bigl\{\bigl\vert
B_1(s_k)\bigr\vert\le a^{-1/\alpha_j} |{\cal
A}_{s_{k-1}} \bigr\}\,da.
\end{eqnarray*}

By Anderson's inequality,
\begin{eqnarray*}&&\P_0 \bigl\{\vert B_{s_k}
\vert\le a^{-1/\alpha_j} |{\cal A}_{s_{k-1}} \bigr\}\\
&&\qquad=\P_0
\bigl\{\bigl\vert B_1(s_{k-1}) +\bigl(B_1(s_k)-B_1(s_{k-1})
\bigr)\bigr\vert\le a^{-1/\alpha_j} |{\cal A}_{s_{k-1}} \bigr\}
\\
&&\qquad\le
\P_0 \bigl\{\bigl\vert B_1(s_k)-B_1(s_{k-1})
\bigr\vert\le a^{-1/\alpha_j} |{\cal A}_{s_{k-1}} \bigr\} =\P_0
\bigl\{\bigl\vert B_1(s_k-s_{k-1})
\bigr\vert^{-\alpha_j}\ge a\bigr\}.
\end{eqnarray*}
So we have
\begin{eqnarray*} \E_0\prod_{k=1}^m
\bigl\vert B_1(s_k)\bigr\vert^{-\alpha_j}&\le&\prod
_{k=1}^m\E_0 \bigl\vert
B_1(s_k-s_{k-1})\bigr\vert^{-\alpha_j}
\\
&=& \bigl\{\E_0 \bigl\vert B_1(1)\bigr\vert^{-\alpha_j} \bigr
\}^m \prod_{k=1}^m(s_k-s_{k-1})^{-\alpha_j},\qquad
j=1,\ldots, d.
\end{eqnarray*}
Here the convention $s_0=0$ is adopted.

Summarizing our computation,
\[
\E_0 \eta_t^m\le m! \Biggl(\prod
_{j=1}^d\E_0 \bigl\vert B_1(1)
\bigr\vert^{-\alpha_j} \Biggr)^m \int_{[0,t]_<^m}\prod
_{k=1}^m(s_k-s_{k-1})^{-\alpha/2}\,ds_1
\cdots ds_m.\vadjust{\goodbreak}
\]

Let $\tau$ be an exponential time with parameter 1 such that
$\tau$ is independent of $B_t$. By Fubini's theorem
%
%
%e5.27 #&#
%
\begin{eqnarray}
\label{a-16} &&\E^\tau\otimes\E_0\eta_\tau^m
\nonumber\\
&&\qquad\le m! \Biggl(\prod_{j=1}^d
\E_0 \bigl\vert B_1(1)\bigr\vert^{-\alpha_j}
\Biggr)^m \nonumber\\
&&\qquad\quad{}\times\int_0^\infty
e^{-t} \Biggl[ \int_{[0,t]_<^m}\prod
_{k=1}^m(s_k-s_{k-1})^{-\alpha/2}\,ds_1
\cdots ds_m \Biggr]\,dt
\\
&&\qquad=m! \Biggl(\prod_{j=1}^d
\E_0 \bigl\vert B_1(1)\bigr\vert^{-\alpha_j}
\Biggr)^m \biggl(\int_0^\infty
t^{-\alpha/2}e^{-t}\,dt \biggr)^m\nonumber \\
&&\qquad=m! \Biggl(\Gamma
\biggl(\frac{2-\alpha}{2} \biggr) \prod_{j=1}^d
\E_0 \bigl\vert B_1(1)\bigr\vert^{-\alpha_j}
\Biggr)^m
\nonumber
\end{eqnarray}
for $m=1,2,\ldots.$

On the other hand, notice that
$ \eta_t\stackrel {d}{=}t^{{(2-\alpha)}/{2}}\eta_1$.
So we have
\[
\E^\tau\otimes\E_0\eta_\tau^m=
\bigl(\E^\tau\tau^{({(2-\alpha)}/{2})m} \bigr) \E_0
\eta_1^m=\Gamma \biggl(1+\frac{2-\alpha}{2}m \biggr)
\E_0\eta_1^m.
\]
Combining this with (\ref{a-16}), by Stirling formula we conclude that
there is a constant $C>0$ such that
\[
\E_0\eta_1^m\le(m!)^{\alpha/2}C^m,\qquad
m=1,2,\ldots.
\]
This implies (\ref{a-15}) with $b<C^{-2/\alpha}$.

We now claim that
%
%
%e5.28 #&#
%
\begin{equation}
\label{a-17} \limsup_{t\to\infty}\frac{1}{ t}\log
\E_0\exp \{\theta\eta _t \} <\infty \qquad\forall\theta>0.
\end{equation}

Indeed, by scaling,
\begin{eqnarray*}\E_0\exp \{\theta\eta_t\} &=&
\E_0\exp \bigl\{\theta t^{{(2-\alpha)}/ {2}}
\eta_1 \bigr\}
\\
&\le&\E_0\exp \bigl\{b\eta_1^{2/\alpha} \bigr\}+
\E_0 \bigl\{\exp \bigl\{\theta t^{{(2-\alpha)} /{2}}
\eta_1 \bigr\}; \eta_1\le\bigl(\theta b^{-1}
\bigr)^{{2}/ {(2-\alpha) }}t^{\alpha/2} \bigr\}
\\
&\le&\E_0\exp \bigl\{b\eta_1^{2/\alpha} \bigr\}+\exp
\bigl\{ \bigl(\theta b^{-1}\bigr)^{{2}/ {(2-\alpha)}} t \bigr
\}. \end{eqnarray*}
Hence, (\ref{a-17}) follows from (\ref{a-15}).

Given $N>0$,
\[
\eta_t\ge\int_0^t
\bigl(K(B_s)\wedge N \bigr)\,ds.
\]
On the other hand, applying Theorem~4.1.6, \cite{Chen} to the bounded,
continuous function $K(x)\wedge N$ gives
\begin{eqnarray*}&&\lim_{t\to\infty}\frac{1}{ t}\log
\E_0\exp \biggl\{\int_0^t
\bigl(K(B_s)\wedge N \bigr)\,ds \biggr\}
\\
&&\qquad=\sup_{g\in{\cal F}_d(\R^d)} \biggl\{\int_{\R^d}
\bigl(K(x)\wedge N \bigr)g^2(x)\,dx -\frac{1}{2}\int
_{\R^d}\bigl\vert\nabla g(x)\bigr\vert^2\,dx \biggr\}.
\end{eqnarray*}
Thus,
\begin{eqnarray*}
&&\sup_{g\in{\cal F}_d(\R^d)} \biggl\{\int_{\R^d} \bigl(K(x)
\wedge N \bigr)g^2(x)\,dx -\frac{1}{2}\int_{\R^d}
\bigl\vert\nabla g(x)\bigr\vert^2\,dx \biggr\}\\
&&\qquad \le\limsup_{t\to\infty}
\frac{1}{ t}\log\E_0\exp \{\theta\eta _t \}.
\end{eqnarray*}
Letting $N\to\infty$ on the left-hand side, by (\ref{a-17}) we have
(\ref{a-13}). \end{pf}

With (\ref{a-13}),
an obvious modification of the argument for (\ref{a-7}) shows that
there is a constant $\widetilde{C}_\alpha>0$ such that the inequality
%
%
%e5.29 #&#
%
\begin{equation}
\label{a-18} \qquad\int_{\R^d} \Biggl(\prod
_{j=1}^d\vert x_j
\vert^{-\alpha_j} \Biggr) f^2(x)\,dx \le\widetilde{C}_\alpha
\|f\|_2^{2-\alpha} \|\nabla f\|_2^{\alpha},\qquad f\in
W^{1,2}\bigl(\R^d\bigr)
\end{equation}
holds. Recall our discussion based on the inequality (\ref{a-1}). Replacing
(\ref{a-1}) by (\ref{a-18}) and copying the same derivation
we obtain a parallel system of inequalities and relations among
variations that are summarized in the following.

First, we have the inequality
%
%
%e5.30 #&#
%
\begin{eqnarray}
\label{a-20}&& \int_{\R^d\times\R^d} \Biggl(\prod
_{j=1}^d\vert x_j-y_j
\vert^{-\alpha_j} \Biggr) f^2(x)f^2(y)\,dx\,dy
\nonumber
\\[-8pt]
\\[-8pt]
\nonumber
&&\qquad\le
\widetilde{C}_\alpha\|f\|_2^{4-\alpha} \|\nabla f
\|_2^{\alpha},\qquad f\in W^{1,2}\bigl(\R^d\bigr).
\end{eqnarray}
Consequently, the best consequence
%
%
%e5.31 #&#
%
\begin{eqnarray}
\label{a-21} &&\tilde{\kappa}(d,\alpha)=\inf \Biggl\{C>0; \int
_{\R^d\times\R^d} \Biggl(\prod_{j=1}^d
\vert x_j-y_j\vert^{-\alpha_j} \Biggr)
f^2(x)f^2(y)\,dx\,dy
\nonumber
\\[-8pt]
\\[-8pt]
\nonumber
&&\hspace*{142pt}\le C\|f\|_2^{4-\alpha} \|\nabla f\|_2^{\alpha}
\ \forall f\in W^{1,2}\bigl(\R^d\bigr) \Biggr\}
\end{eqnarray}
is finite.\vadjust{\goodbreak}

Second, the quantities defined through the variations
%
%
%e5.32 #&#
%
\begin{eqnarray}
\label{a-22} &&\widetilde{M}_{d,\alpha}(\theta)\nonumber\\
&&\qquad=\sup_{g\in{\cal F}_d(\R^d)}
\Biggl\{\theta \Biggl( \int_{\R^d\times\R^d} \Biggl(\prod
_{j=1}^d\vert x_j-y_j
\vert ^{-\alpha
_j} \Biggr) g^2(x)g^2(y)\,dx\,dy
\Biggr)^{1/2}
\nonumber
\\[-8pt]
\\[-8pt]
\nonumber&&\hspace*{211pt}{}-\frac{1}{2}\int_{\R^d}\bigl\vert\nabla g(x)
\bigr\vert^2\,dx \Biggr\}, \\
\eqntext{\theta>0,}\\
\label{a-23} &&\tilde{\sigma}(d,\alpha)
\nonumber
\\[-8pt]
\\[-8pt]
\nonumber
\qquad&&\qquad=\sup_{f\in{\cal G}_d(\R^d)} \Biggl\{ \int
_{\R^d\times\R^d} \Biggl(\prod_{j=1}^d
\vert x_j-y_j\vert ^{-\alpha
_j} \Biggr)
f^2(x)f^2(y)\,dx\,dy \Biggr\}^{1/2}
\end{eqnarray}
are finite

Third, these variations are co-related according to the following lemma.

%
%le5.6 #&#
\begin{lemmas}\label{a-24} Under $0\le\alpha_j<1$ ($j=1,\ldots, d$
and $\alpha_1+\cdots+\alpha_d <2$,
%
%
%e5.34 #&#
%
\begin{eqnarray}
\label{a-25} \widetilde{M}_{d,\alpha}(\theta)&=&\frac{4-\alpha}{4} \biggl(
\frac{\alpha}{2} \biggr)^{{\alpha} /{(4-\alpha)}} \tilde{\kappa}(d,
\alpha)^{{2} /{(4-\alpha)}}\theta^{{4} /{(4-\alpha)}},
\\
\label{a-26} \tilde{\sigma}(d,\alpha)&=& \biggl(\frac{4-\alpha}{4}
\biggr)^{{(4-\alpha) }/ {4}} \biggl(\frac{\alpha}{2}
\biggr)^{\alpha/4} \tilde{\kappa}(d,\alpha)^{1/2}.
\end{eqnarray}
\end{lemmas}

The next lemma is related to Theorem~\ref{intro-46}.

%
%le5.7 #&#
\begin{lemmas}\label{a-27}
%
%
%e5.36 #&#
%
\begin{eqnarray}
\label{a-28} &&\sup_{g\in{\cal F}_1(\R)} \biggl\{\theta \biggl(\int
_{-\infty
}^\infty g^4(x)\,dx
\biggr)^{1/2} -\frac{1}{2}\int_{-\infty}^\infty
\bigl\vert f'(x)\bigr\vert^2 \,dx \biggr\}
\nonumber
\\[-8pt]
\\[-8pt]
\nonumber
&&\qquad=\frac{1}{2}
\biggl(\frac{3}{4} \biggr)^{2/3}\theta^{4/3}\qquad (\theta>0),
\\
\label{a-29} &&\sup_{g\in{\cal G}_1(\R)}\int_{-\infty}^\infty
g^4(x)\,dx =\frac{3}{4} \biggl(\frac{1}{2}
\biggr)^{3/2}.
\end{eqnarray}
\end{lemmas}

\begin{pf} The identity (\ref{a-28}) is given in Theorem C.4, page
307, \cite{Chen}.
This theorem also claims the Sobolev inequality
\[
\|f\|_4\le3^{-1/8}\|f\|_2^{3/4}
\bigl\|f'\bigr\|_2^{3/4},\qquad f\in W^{1,2}\bigl(
\R^d\bigr)
\]
with $3^{-1/8}$ as the best constant. A natural
modification of the proof for (\ref{a-8})
leads to (\ref{a-29}).
\end{pf}
\end{appendix}

\section*{Acknowledgments} The author
would like to thank
two anonymous referees who read the first version
of this paper for their insightful remarks, critical comments and
excellent suggestions.

%
% imsref loaded by akundreckaite, 2013-08-26 10:35:21


%suskaldyti doi

\printaddresses


\begin{thebibliography}{26}
% BibTex style file: ims.bst, 2013-01-28
% Default style options (sort=0,type=number).
% Used options (sort=1,type=number).

%b1 #&#
\bibitem{ACQ}
\begin{barticle}[mr]
\bauthor{\bsnm{Amir},~\bfnm{Gideon}\binits{G.}},
  \bauthor{\bsnm{Corwin},~\bfnm{Ivan}\binits{I.}} \AND
  \bauthor{\bsnm{Quastel},~\bfnm{Jeremy}\binits{J.}}
(\byear{2011}).
\btitle{Probability distribution of the free energy of the continuum directed
  random polymer in {$1+1$} dimensions}.
\bjournal{Comm. Pure Appl. Math.}
\bvolume{64}
\bpages{466--537}.
\bid{doi={10.1002/cpa.20347}, issn={0010-3640}, mr={2796514}}
\bptok{imsref}%
\end{barticle}
\endbibitem

%b2 #&#
\bibitem{BQS}
\begin{barticle}[mr]
\bauthor{\bsnm{Bal{\'a}zs},~\bfnm{M.}\binits{M.}},
  \bauthor{\bsnm{Quastel},~\bfnm{J.}\binits{J.}} \AND
  \bauthor{\bsnm{Sepp{\"a}l{\"a}inen},~\bfnm{T.}\binits{T.}}
(\byear{2011}).
\btitle{Fluctuation exponent of the {KPZ}/stochastic {B}urgers equation}.
\bjournal{J. Amer. Math. Soc.}
\bvolume{24}
\bpages{683--708}.
\bid{doi={10.1090/S0894-0347-2011-00692-9}, issn={0894-0347}, mr={2784327}}
\bptok{imsref}%
\end{barticle}
\endbibitem

%b3 #&#
\bibitem{BCR}
\begin{barticle}[mr]
\bauthor{\bsnm{Bass},~\bfnm{Richard}\binits{R.}},
  \bauthor{\bsnm{Chen},~\bfnm{Xia}\binits{X.}} \AND
  \bauthor{\bsnm{Rosen},~\bfnm{Jay}\binits{J.}}
(\byear{2009}).
\btitle{Large deviations for {R}iesz potentials of additive processes}.
\bjournal{Ann. Inst. Henri Poincar\'e Probab. Stat.}
\bvolume{45}
\bpages{626--666}.
\bid{doi={10.1214/08-AIHP181}, issn={0246-0203}, mr={2548497}}
\bptok{imsref}%
\end{barticle}
\endbibitem

%b4 #&#
\bibitem{BK1}
\begin{barticle}[mr]
\bauthor{\bsnm{Biskup},~\bfnm{Marek}\binits{M.}} \AND
  \bauthor{\bsnm{K{\"o}nig},~\bfnm{Wolfgang}\binits{W.}}
(\byear{2001}).
\btitle{Long-time tails in the parabolic {A}nderson model with bounded
  potential}.
\bjournal{Ann. Probab.}
\bvolume{29}
\bpages{636--682}.
\bid{doi={10.1214/aop/1008956688}, issn={0091-1798}, mr={1849173}}
\bptok{imsref}%
\end{barticle}
\endbibitem

%b5 #&#
\bibitem{CM-1}
\begin{barticle}[mr]
\bauthor{\bsnm{Carmona},~\bfnm{R.~A.}\binits{R.~A.}} \AND
  \bauthor{\bsnm{Molchanov},~\bfnm{S.~A.}\binits{S.~A.}}
(\byear{1995}).
\btitle{Stationary parabolic {A}nderson model and intermittency}.
\bjournal{Probab. Theory Related Fields}
\bvolume{102}
\bpages{433--453}.
\bid{doi={10.1007/BF01198845}, issn={0178-8051}, mr={1346261}}
\bptok{imsref}%
\end{barticle}
\endbibitem

%b6 #&#
\bibitem{CV}
\begin{barticle}[mr]
\bauthor{\bsnm{Carmona},~\bfnm{Ren{\'e}~A.}\binits{R.~A.}} \AND
  \bauthor{\bsnm{Viens},~\bfnm{Frederi~G.}\binits{F.~G.}}
(\byear{1998}).
\btitle{Almost-sure exponential behavior of a stochastic {A}nderson model with
  continuous space parameter}.
\bjournal{Stochastics Stochastics Rep.}
\bvolume{62}
\bpages{251--273}.
\bid{issn={1045-1129}, mr={1615092}}
\bptok{imsref}%
\end{barticle}
\endbibitem

%b7 #&#
\bibitem{Chen}
\begin{bbook}[mr]
\bauthor{\bsnm{Chen},~\bfnm{Xia}\binits{X.}}
(\byear{2010}).
\btitle{Random Walk Intersections: Large Deviations and Related Topics}.
\bseries{Mathematical Surveys and Monographs}
\bvolume{157}.
\bpublisher{Amer. Math. Soc.}, \blocation{Providence, RI}.
\bid{mr={2584458}}
\bptnote{check year}%
\bptok{imsref}%
\end{bbook}
\endbibitem

%b8 #&#
\bibitem{C}
\begin{barticle}[mr]
\bauthor{\bsnm{Chen},~\bfnm{Xia}\binits{X.}}
(\byear{2012}).
\btitle{Quenched asymptotics for {B}rownian motion of renormalized {P}oisson
  potential and for the related parabolic {A}nderson models}.
\bjournal{Ann. Probab.}
\bvolume{40}
\bpages{1436--1482}.
\bid{doi={10.1214/11-AOP655}, issn={0091-1798}, mr={2978130}}
\bptok{imsref}%
\end{barticle}
\endbibitem

%b9 #&#
\bibitem{CHSX}
\begin{bmisc}[auto:STB|2013/06/05|13:45:01]
\bauthor{\bsnm{Chen},~\bfnm{X.}\binits{X.}},
  \bauthor{\bsnm{Hu},~\bfnm{Y.~Z.}\binits{Y.~Z.}},
  \bauthor{\bsnm{Song},~\bfnm{J.}\binits{J.}} \AND
  \bauthor{\bsnm{Xing},~\bfnm{F.}\binits{F.}}
  (\byear{2014}).
\bhowpublished{Exponential asymptotics for time-space Hamiltonians.
\textit{Ann. Inst. Henri Poincar\'e Probab. Stat.} To appear}.
\bptok{imsref}%
\end{bmisc}
\endbibitem

%b10 #&#
\bibitem{CLR}
\begin{barticle}[mr]
\bauthor{\bsnm{Chen},~\bfnm{Xia}\binits{X.}},
  \bauthor{\bsnm{Li},~\bfnm{Wenbo~V.}\binits{W.~V.}} \AND
  \bauthor{\bsnm{Rosen},~\bfnm{Jay}\binits{J.}}
(\byear{2005}).
\btitle{Large deviations for local times of stable processes and stable random
  walks in 1 dimension}.
\bjournal{Electron. J. Probab.}
\bvolume{10}
\bpages{577--608}.
\bid{doi={10.1214/EJP.v10-260}, issn={1083-6489}, mr={2147318}}
\bptok{imsref}%
\end{barticle}
\endbibitem

%b11 #&#
\bibitem{C-R}
\begin{barticle}[mr]
\bauthor{\bsnm{Chen},~\bfnm{Xia}\binits{X.}} \AND
  \bauthor{\bsnm{Rosen},~\bfnm{Jay}\binits{J.}}
(\byear{2010}).
\btitle{Large deviations and renormalization for {R}iesz potentials of stable
  intersection measures}.
\bjournal{Stochastic Process. Appl.}
\bvolume{120}
\bpages{1837--1878}.
\bid{doi={10.1016/j.spa.2010.05.006}, issn={0304-4149}, mr={2673977}}
\bptok{imsref}%
\end{barticle}
\endbibitem

%b12 #&#
\bibitem{CJKS}
\begin{barticle}[mr]
\bauthor{\bsnm{Conus},~\bfnm{Daniel}\binits{D.}},
  \bauthor{\bsnm{Joseph},~\bfnm{Mathew}\binits{M.}},
  \bauthor{\bsnm{Khoshnevisan},~\bfnm{Davar}\binits{D.}} \AND
  \bauthor{\bsnm{Shiu},~\bfnm{Shang-Yuan}\binits{S.-Y.}}
(\byear{2013}).
\btitle{On the chaotic character of the stochastic heat equation, {II}}.
\bjournal{Probab. Theory Related Fields}
\bvolume{156}
\bpages{483--533}.
\bid{doi={10.1007/s00440-012-0434-3}, issn={0178-8051}, mr={3078278}}
\bptnote{check year}%
\bptok{imsref}%
\end{barticle}
\endbibitem

%b13 #&#
\bibitem{G-K}
\begin{barticle}[mr]
\bauthor{\bsnm{G{\"a}rtner},~\bfnm{J{\"u}rgen}\binits{J.}} \AND
  \bauthor{\bsnm{K{\"o}nig},~\bfnm{Wolfgang}\binits{W.}}
(\byear{2000}).
\btitle{Moment asymptotics for the continuous parabolic {A}nderson model}.
\bjournal{Ann. Appl. Probab.}
\bvolume{10}
\bpages{192--217}.
\bid{doi={10.1214/aoap/1019737669}, issn={1050-5164}, mr={1765208}}
\bptok{imsref}%
\end{barticle}
\endbibitem

%b14 #&#
\bibitem{G-K-M}
\begin{barticle}[mr]
\bauthor{\bsnm{G{\"a}rtner},~\bfnm{J.}\binits{J.}},
  \bauthor{\bsnm{K{\"o}nig},~\bfnm{W.}\binits{W.}} \AND
  \bauthor{\bsnm{Molchanov},~\bfnm{S.~A.}\binits{S.~A.}}
(\byear{2000}).
\btitle{Almost sure asymptotics for the continuous parabolic {A}nderson model}.
\bjournal{Probab. Theory Related Fields}
\bvolume{118}
\bpages{547--573}.
\bid{doi={10.1007/PL00008754}, issn={0178-8051}, mr={1808375}}
\bptok{imsref}%
\end{barticle}
\endbibitem

%b15 #&#
\bibitem{G-M}
\begin{barticle}[mr]
\bauthor{\bsnm{G{\"a}rtner},~\bfnm{J.}\binits{J.}} \AND
  \bauthor{\bsnm{Molchanov},~\bfnm{S.~A.}\binits{S.~A.}}
(\byear{1990}).
\btitle{Parabolic problems for the {A}nderson model. {I}. {I}ntermittency and
  related topics}.
\bjournal{Comm. Math. Phys.}
\bvolume{132}
\bpages{613--655}.
\bid{issn={0010-3616}, mr={1069840}}
\bptok{imsref}%
\end{barticle}
\endbibitem

%b16 #&#
\bibitem{G-M-1}
\begin{barticle}[mr]
\bauthor{\bsnm{G{\"a}rtner},~\bfnm{J.}\binits{J.}} \AND
  \bauthor{\bsnm{Molchanov},~\bfnm{S.~A.}\binits{S.~A.}}
(\byear{1998}).
\btitle{Parabolic problems for the {A}nderson model. {II}. {S}econd-order
  asymptotics and structure of high peaks}.
\bjournal{Probab. Theory Related Fields}
\bvolume{111}
\bpages{17--55}.
\bid{doi={10.1007/s004400050161}, issn={0178-8051}, mr={1626766}}
\bptok{imsref}%
\end{barticle}
\endbibitem

%b17 #&#
\bibitem{GV}
\begin{bbook}[auto:STB|2013/06/05|13:45:01]
\bauthor{\bsnm{Guelfand},~\bfnm{I.~M.}\binits{I.~M.}} \AND
  \bauthor{\bsnm{Vilenkin},~\bfnm{G.}\binits{G.}}
(\byear{1964}).
\btitle{Generalized Functions}.
\bpublisher{Academic Press}, \blocation{New York}.
\bptok{imsref}%
\end{bbook}
\endbibitem

%b18 #&#
\bibitem{H}
\begin{barticle}[mr]
\bauthor{\bsnm{Hairer},~\bfnm{M.}\binits{M.}}
(\byear{2013}).
\btitle{Solving the KPZ equation}.
\bjournal{Ann. Math.}
\bvolume{178}
\bpages{559--664}.
\bid{mr={3071506}}
\bptok{imsref}%
\end{barticle}
\endbibitem

%b19 #&#
\bibitem{HNS}
\begin{barticle}[mr]
\bauthor{\bsnm{Hu},~\bfnm{Yaozhong}\binits{Y.}},
  \bauthor{\bsnm{Nualart},~\bfnm{David}\binits{D.}} \AND
  \bauthor{\bsnm{Song},~\bfnm{Jian}\binits{J.}}
(\byear{2011}).
\btitle{Feynman--{K}ac formula for heat equation driven by fractional white
  noise}.
\bjournal{Ann. Probab.}
\bvolume{39}
\bpages{291--326}.
\bid{doi={10.1214/10-AOP547}, issn={0091-1798}, mr={2778803}}
\bptok{imsref}%
\end{barticle}
\endbibitem

%b20 #&#
\bibitem{karda-1}
\begin{barticle}[auto:STB|2013/06/05|13:45:01]
\bauthor{\bsnm{Karda},~\bfnm{M.}\binits{M.}},
  \bauthor{\bsnm{Parisi},~\bfnm{G.}\binits{G.}} \AND
  \bauthor{\bsnm{Zhang},~\bfnm{Y.~C.}\binits{Y.~C.}}
(\byear{1986}).
\btitle{Dynamic scaling of growing interface}.
\bjournal{Phys. Rev. Lett.}
\bvolume{56}
\bpages{889--892}.
\bptok{imsref}%
\end{barticle}
\endbibitem

%b21 #&#
\bibitem{karda-2}
\begin{barticle}[auto:STB|2013/06/05|13:45:01]
\bauthor{\bsnm{Karda},~\bfnm{M.}\binits{M.}} \AND
  \bauthor{\bsnm{Zhang},~\bfnm{Y.~C.}\binits{Y.~C.}}
(\byear{1987}).
\btitle{Scaling of directed polymers in random media}.
\bjournal{Phys. Rev. Lett.}
\bvolume{58}
\bpages{2087--2090}.
\bptok{imsref}%
\end{barticle}
\endbibitem

%b22 #&#
\bibitem{MR}
\begin{bbook}[mr]
\bauthor{\bsnm{Marcus},~\bfnm{Michael~B.}\binits{M.~B.}} \AND
  \bauthor{\bsnm{Rosen},~\bfnm{Jay}\binits{J.}}
(\byear{2006}).
\btitle{Markov Processes, {G}aussian Processes, and Local Times}.
\bseries{Cambridge Studies in Advanced Mathematics}
\bvolume{100}.
\bpublisher{Cambridge Univ. Press}, \blocation{Cambridge}.
\bid{doi={10.1017/CBO9780511617997}, mr={2250510}}
\bptok{imsref}%
\end{bbook}
\endbibitem

%b23 #&#
\bibitem{Slepian}
\begin{barticle}[mr]
\bauthor{\bsnm{Slepian},~\bfnm{David}\binits{D.}}
(\byear{1962}).
\btitle{The one-sided barrier problem for {G}aussian noise}.
\bjournal{Bell System Tech. J.}
\bvolume{41}
\bpages{463--501}.
\bid{issn={0005-8580}, mr={0133183}}
\bptok{imsref}%
\end{barticle}
\endbibitem

%b24 #&#
\bibitem{Sznitman}
\begin{bbook}[mr]
\bauthor{\bsnm{Sznitman},~\bfnm{Alain-Sol}\binits{A.-S.}}
(\byear{1998}).
\btitle{Brownian Motion, Obstacles and Random Media}.
\bpublisher{Springer}, \blocation{Berlin}.
\bid{mr={1717054}}
\bptok{imsref}%
\end{bbook}
\endbibitem

%b25 #&#
\bibitem{VZ}
\begin{barticle}[mr]
\bauthor{\bsnm{Viens},~\bfnm{Frederi~G.}\binits{F.~G.}} \AND
  \bauthor{\bsnm{Zhang},~\bfnm{Tao}\binits{T.}}
(\byear{2008}).
\btitle{Almost sure exponential behavior of a directed polymer in a fractional
  {B}rownian environment}.
\bjournal{J. Funct. Anal.}
\bvolume{255}
\bpages{2810--2860}.
\bid{doi={10.1016/j.jfa.2008.06.020}, issn={0022-1236}, mr={2464192}}
\bptok{imsref}%
\end{barticle}
\endbibitem

\end{thebibliography}
\end{document}